\documentclass[article,preprint,a4paper,10pt,times]{article}
\usepackage{amsthm,amscd,amsfonts,amssymb,graphicx}
\usepackage{mathtools,mathrsfs}
\usepackage[bookmarksnumbered,colorlinks,plainpages]{hyperref}

\usepackage{geometry}
\DeclareMathAlphabet\EuScript{U}{eus}{m}{n}
\DeclareMathAlphabet\EuScriptBold{U}{eus}{b}{n}
\DeclareMathAlphabet\EuRom{U}{eur}{m}{n}
\DeclareMathAlphabet\EuRomBold{U}{eur}{b}{n}
\DeclareMathAlphabet\EuFrak{U}{euf}{m}{n}
\DeclareMathAlphabet\EuFrakBold{U}{euf}{b}{n}
\DeclareMathAlphabet\EuScriptItal{U}{eus}{i}{i}

\newcommand\mathcalb{\EuScriptBold}

\let\CMcal=\mathcal

\newcommand\TS{\textstyle}
\newcommand\IS{\scriptstyle}
\newcommand\IIS{\scriptscriptstyle}
\newcommand\sstyle{\scriptstyle}
\newcommand\ssstyle{\scriptscriptstyle}

\theoremstyle{plain}
\newtheorem{theorem}{Theorem}[section]
\newtheorem{lemma}[theorem]{Lemma}
\newtheorem{corollary}[theorem]{Corollary}
\newtheorem{proposition}[theorem]{Proposition}

\theoremstyle{definition}
\newtheorem{definition}[theorem]{Definition}
\newtheorem{example}[theorem]{Example}

\theoremstyle{remark}
\newtheorem{remark}{\rm Remark}

\newcommand\bbproo{\noindent{\it Proof:\/} }
\newcommand\eeproo{$ $\hfill $\square$}
\newcommand\bbdef{\begin{definition}\sl }
\newcommand\eedef{\end{definition}\rm }
\newcommand\bbt{\begin{theorem}}
\newcommand\eet{\end{theorem}}
\newcommand\bbc{\begin{corollary}}
\newcommand\eec{\end{corollary}}
\newcommand\bbr{\begin{remark}\rm}
 \newcommand\eer{\end{remark}}
 \newcommand\bbex{\begin{example}\rm}
 \newcommand\eeex{\end{example}}
 \newcommand\bbbt{\begin{theorem}}
 \newcommand\eeet{\end{theorem}}
 \newcommand\bbbl{\begin{lemma}}
 \newcommand\eeel{\end{lemma}}
 \newcommand\bbbc{\begin{corollary}}
 \newcommand\eeec{\end{corollary}}
 \newcommand\bbprop{\begin{proposition}}
 \newcommand\eeprop{\end{proposition}}
 \newcommand\bbl{\begin{lemma}}
 \newcommand\eel{\end{lemma}}
\newcommand\bbn{\begin{enumerate}}
\newcommand\een{\end{enumerate}}

\newcommand\bbe{\begin{equation}}   
\newcommand\eeE{\end{equation}}

\newcommand\bbbe{\begin{equation}}
\newcommand\eeee{\end{equation}}
\newcommand\bbbee{\begin{eqnarray}} 
\newcommand\eeeee{\end{eqnarray}}
\newcommand\BBBB{\begin{eqnarray*}}
\newcommand\EEEE{\end{eqnarray*}}
\newcommand\bba{\begin{eqnarray}}
\newcommand\eea{\end{eqnarray}}
\newcommand\bbz{\begin{eqnarray*}}
\newcommand\eez{\end{eqnarray*}}

\newcommand\bbal{\begin{align}}               
\newcommand\eeal{\end{align}}

\newcommand\vsplind{\vspace{-2mm}\newline}
\newcommand\vsplindp{\vspace{-2.5mm}\newline}
\newcommand\vsplint{\vspace{-3mm}\newline}
\newcommand\vsplintp{\vspace{-3.5mm}\newline}
\newcommand\vsplinc{\vspace{-4mm}\newline}
\newcommand\vsplincp{\vspace{-4.5mm}\newline}
\newcommand\vsplinp{\vspace{-5mm}\newline}

\newcommand\vsplins{\vspace{-7mm}\newline}

\newcommand\vsplinix{\vspace{-9mm}\newline}
\newcommand\vsplinx{\vspace{-1cm}\newline}
\newcommand\vsplinxi{\vspace{-11mm}\newline}

\newcommand\bbalvsp[1]{\vspace{-#1mm}\begin{align}}
\newcommand\bbalzvsp[1]{\vspace{-#1mm}\begin{align*}}
\newcommand\vspbbal[1]{\vspace{-#1mm}\begin{align}}
\newcommand\vspbbalj{\vspace{-1mm}\begin{align}}
\newcommand\vspbbaljp{\vspace{-1.5mm}\begin{align}}
\newcommand\vspbbald{\vspace{-2mm}\begin{align}}
\newcommand\vspbbaldj{\vspace{-2.1mm}\begin{align}}
\newcommand\vspbbaldd{\vspace{-2.2mm}\begin{align}}
\newcommand\vspbbaldt{\vspace{-2.3mm}\begin{align}}
\newcommand\vspbbaldc{\vspace{-2.4mm}\begin{align}}
\newcommand\vspbbaldp{\vspace{-2.5mm}\begin{align}}
\newcommand\vspbbalt{\vspace{-3mm}\begin{align}}
\newcommand\vspbbaltp{\vspace{-3.5mm}\begin{align}}
\newcommand\vspbbalc{\vspace{-4mm}\begin{align}}
\newcommand\vspbbalcp{\vspace{-4.5mm}\begin{align}}
\newcommand\vspbbalp{\vspace{-5mm}\begin{align}}
\newcommand\vspbbalpp{\vspace{-5.5mm}\begin{align}}
\newcommand\vspbbalvi{\vspace{-6mm}\begin{align}}
\newcommand\vspbbalvip{\vspace{-6.5mm}\begin{align}}
\newcommand\vspbbals{\vspace{-7mm}\begin{align}}
\newcommand\vspbbalsp{\vspace{-7.5mm}\begin{align}}
\newcommand\vspbbalo{\vspace{-8mm}\begin{align}}
\newcommand\vspbbalop{\vspace{-8.5mm}\begin{align}}
\newcommand\vspbbalz[1]{\vspace{#1mm}\begin{align*}}
\newcommand\vspbbalzj{\vspace{-1mm}\begin{align*}}
\newcommand\vspbbalzjp{\vspace{-1.5mm}\begin{align*}}
\newcommand\vspbbalzd{\vspace{-2mm}\begin{align*}}
\newcommand\vspbbalzdp{\vspace{-2.5mm}\begin{align*}}
\newcommand\vspbbalzt{\vspace{-3mm}\begin{align*}}
\newcommand\vspbbalztp{\vspace{-3.5mm}\begin{align*}}
\newcommand\vspbbalzc{\vspace{-4mm}\begin{align*}}
\newcommand\vspbbalzcp{\vspace{-4.5mm}\begin{align*}}
\newcommand\vspbbalzp{\vspace{-5mm}\begin{align*}}
\newcommand\vspbbalzpp{\vspace{-5.5mm}\begin{align*}}
\newcommand\vspbbalzvi{\vspace{-6mm}\begin{align*}}
\newcommand\vspbbalzvip{\vspace{-6.5mm}\begin{align*}}
\newcommand\vspbbalzs{\vspace{-7mm}\begin{align*}}
\newcommand\vspbbalzsp{\vspace{-7.5mm}\begin{align*}}

\newcommand\bbalz[1]{\vspace{-#1mm}\begin{align*}}
\newcommand\bbeez[1]{\vspace{-#1mm}\begin{equation*}}
\newcommand\bbee[1]{\vspace{-#1mm}\begin{equation}}
\newcommand\bbalu[1]{\vspace{-#1mm}\begin{align}}
\newcommand\bbalj{\vspace{-1mm}\begin{align}}
\newcommand\bbald{\vspace{-2mm}\begin{align}}
\newcommand\bbalt{\vspace{-3mm}\begin{align}}
\newcommand\bbalc{\vspace{-4mm}\begin{align}}
\newcommand\bbalzo{\begin{align*}}
\newcommand\bbalzj{\vspace{-1mm}\begin{align*}}
\newcommand\bbalzd{\vspace{-2mm}\begin{align*}}
\newcommand\bbalzt{\vspace{-3mm}\begin{align*}}
\newcommand\bbalzc{\vspace{-4mm}\begin{align*}}

\newcommand\eR{{\mathbb R}}

\newcommand\Rp{{\eR_{\negj\IIS+}}}

\newcommand\Rpz{\eR_{\negj\IIS+}^{\negj*}}

\newcommand\eRp{{[0,\!+\iii\negd)}}

\newcommand\eRpo{{(\negj 0,\negd+\iii\negd)}}
\newcommand\eRpz{{(\negj 0,\!+\iii\negd)}}

\newcommand\Ce{{\mathbb C}}

\newcommand\eN{{\mathbb N}}

\newcommand\en{{\mathbb N}}

\newcommand\OJ{{\negj[\negj0\negj,\negd1\negj]}}

\newcommand\OJo{{(\negj0\negj,\negd1\negj)}}

\newcommand\bbSOJoe{\bbS\odu{0,\!1}}
\newcommand\bbSOJo{\bbS\oduu{0,\!1}}

\providecommand\qed\square
\newcommand\alf\alpha
\newcommand\vei\varepsilon
\newcommand\vep\varepsilon
\newcommand\vfi\varphi

\newcommand\les\leqslant
\newcommand\ges\geqslant
\newcommand\iii\infty

\newcommand\ttt\rightarrow
\newcommand\ifff\Leftrightarrow
\newcommand\mto\mapsto
\newcommand\sset\subset
\newcommand\smin\setminus
\newcommand\rest\restriction

\newcommand\ssearrow{{\pozd\ssstyle\searrow\pozj}}

\newcommand\sssearrow{{\pozd\ssstyle\searrow\pozj}}

\newcommand\dir{{\IIS\triangleright}}
\newcommand\ISj{{\negj_{\IS1}}}
\newcommand\ISd{{\negj_{\IS2}}}

\newcommand\x{{\times}}

\newcommand\xuu{\negd\times\negd}

\newcommand\oxuu{\negd\otimes\negd}
\newcommand\oxtj{\!\otimes\negj}

\newcommand\oxdt{\negd\otimes\!}
\newcommand\oxtd{\!\otimes\negd}

\newcommand\omin\ominus

\newcommand\mmm{\textrm{\textrm{-}}}

\newcommand\mmi\mmm

\newcommand\minu{\negj-\negj}
\newcommand\minjd{\negj-\negd}
\newcommand\mindj{\negd-\negj}

\newcommand\minuu{\negd-\negd}

\newcommand\mintj{\!-\negj}

\newcommand\mindt{\negd-\!}

\newcommand\mintt{\!-\!}

\newcommand\minpj{\negp-\negj}
\newcommand\mindp{\negd-\negp}

\newcommand\minct{\negc-\!}
\newcommand\minpd{\negp-\negd}

\newcommand\plu{\negj+\negj}
\newcommand\plujj{\negj+\negj}

\newcommand\pludj{\negd+\negj}
\newcommand\plujt{\negj+\!}
\newcommand\pluu{\negd+\negd}

\newcommand\pludt{\negd+\!}

\newcommand\pludc{\negd+\negc}

\newcommand\pludp{\negd+\negp}
\newcommand\plutc{\!+\negc}

\newcommand\plutp{\!+\negp}

\newcommand\plusdj{\negd+\negj}

\newcommand\plusdc{\negd+\negc}

\newcommand\plusdp{\negd+\negp}

\newcommand\cupu{\negj\cup\negj}

\newcommand\captj{\!\cap\negj}

\newcommand\tr{\textrm{\rm tr}}

\newcommand\kor[1]{\negj\sqrt{\negj#1}}

\newcommand\TSsqrt[1]{{\TS\negj\sqrt{#1\negj}}}

\newcommand\ufrac[2]{{\negj\frac{#1}{#2}}}
\newcommand\uufrac[2]{{\negd\frac{#1}{#2}}}
\newcommand\utfrac[2]{{\!\frac{#1}{#2}}}
\newcommand\ucfrac[2]{{\negc\frac{#1}{#2}}}

\newcommand\uvifrac[2]{{\!\!\frac{#1}{#2}}}
\newcommand\usfrac[2]{{\negs\frac{#1}{#2}}}

\newcommand\kroz[1]{\negj/\negj#1}

\newcommand\eqdu{\negd\stackrel{\textnormal{\tiny def}}{=}\negd}

\newcommand\eqdto{\!\stackrel{\textnormal{\tiny def}}{=}}
\newcommand\eqdtj{\negc\stackrel{\textnormal{\tiny def}}{=}\negd}

\makeatletter
\newcommand\dod{\negj\mathrel{\rlap{\raisebox{0.3ex} {$\m@th\cdot$}}\raisebox{-0.3ex}{$\m@th\cdot$}}=\negj}

\newcommand\dodu{\negd\mathrel{\rlap{\raisebox{0.3ex}{$\m@th\cdot$}}\raisebox{-0.3ex}{$\m@th\cdot$}}=\negd}

\newcommand\eqsimm{\mathrel{\rlap{\raisebox{0.38ex}{$\m@th\sim$}}\raisebox{-0.05ex}{$\m@th-$}\llap{\raisebox{-0.33ex}{$\m@th-$}} }}

\newcommand\eqoj{=\negj}
\newcommand\eqod{=\negd}
\newcommand\equ{\negj=\negj}

\newcommand\eqdj{\negd=\negj}

\newcommand\eqjt{=\!}
\newcommand\equu{\negd=\negd}

\newcommand\eqdt{\negd=\!}

\newcommand\eqcj{\negc=\negj}

\newcommand\eqdc{\negd=\negc}

\newcommand\equt{\!=\!}

\newcommand\eqpj{\negp=\negj}

\newcommand\eqdp{\negd=\negp}

\newcommand\eqpd{\negp=\negd}

\newcommand\equc{\negc=\negc}

\newcommand\equp{\negp=\negp}
\newcommand\eqsd{\negs=\negd}

\newcommand\eqxd{\negx=\negd}

\newcommand\lesoj{\les\negj}             
\newcommand\lesod{\les\negd}
\newcommand\lesu{\negj\les\negj}

\newcommand\lesot{\les\!}

\newcommand\lesdj{\negd\les\negj}
\newcommand\lesoc{\les\negc}
\newcommand\lestj{\!\les\negj}

\newcommand\lescj{\negc\les\negj}
\newcommand\lesuu{\negd\les\negd}

\newcommand\lespj{\negp\les\negj}
\newcommand\lesut{\!\les\!}

\newcommand\lespd{\negp\les\negd}

\newcommand\lespt{\negp\les\!}
\newcommand\lesst{\negs\les\!}

\newcommand\lessd{\negs\les\negd}

\newcommand\lesup{\negp\les\negp}

\newcommand\lesxd{\negx\les\negd}

\newcommand\gesu{\negj\ges\negj}
\newcommand\gesuu{\negd\ges\negd}

\newcommand\lsuu{\negd<\negd}

\newcommand\inu{\negj\in\negj}              

\newcommand\inuu{\negd\in\negd}

\newcommand\intj{\!\in\negj}

\newcommand\inN{\in\mathbb N}

\newcommand\iip{{+\negj\infty}}

\newcommand\GCD{\Gamma_{\ssstyle\negc C\!,\negj D}}

\newcommand\GCzCI[1]{\Gamma_{\ssstyle\negp C^*\nego,\negj C}^{#1}\negd(\negj I\negj)}
\newcommand\GCCzI[1]{\Gamma_{\ssstyle\negc C\negd,\negj C^*\negj}^{#1}\!(\negj I\negj)}
\newcommand\GDzDI[1]{\Gamma_{\ssstyle\negs D^{\negj*}\nego,\negj D}^{#1}\negd(\negj I\negj)}
\newcommand\GDDzI[1]{\Gamma_{\ssstyle\negc D\negd,D^{\negj*}\negj}^{#1}\!(\negj I\negj)}

\newcommand\calAC{A_{\negj_C}}
\newcommand\calACz{A_{\negj_{C^{\negj*}}}\negj}
\newcommand\calAD{A_{\negj_D}}
\newcommand\calADz{A_{\negj_{D^{\negj*}}}\negj}

\newcommand\mj{{\minjd1}}

\newcommand\ujd{{\negd 1\negd/\negj2}}
\newcommand\ejd{{1\negd/\negj2}}
\newcommand\uujd{{\! 1\negd/\negj2}}

\newcommand\utjd{{\negc1\negd/\negj2}}

\newcommand\upjd{{\negvi1\negd/\negj2}}

\newcommand\ucfjd{{\negc\frac12}}

\newcommand\nti{{n\ttt\iii}}             

\newcommand\Anuz{A_n^{\negd*}}
\newcommand\Anff{A_n^\fff}

\newcommand\Bn{B_n^\fff}

\newcommand\Bnuz{B_n^*}

\newcommand\Cn{C_n^\fff}

\newcommand\Dn{D_n^\fff}

\newcommand\AtXBt{A_t^\fff XB_t^\fff}

\newcommand\AnXBn{A_n^\fff\negj X\negj B_{\negj n}^\fff}

\newcommand\CnXDn{C_n^\fff\negj X\negj D_{\negd n}^\fff}

\newcommand\CnzCn{C_n^{\negj*}\negj C_n^\fff}
\newcommand\CnCnz{C_n^\fff\negj C_n^*}
\newcommand\DnzDn{D_n^{\negj*}\negj D_n^\fff}
\newcommand\DnDnz{D_n^\fff\negj D_n^*}

\newcommand\AnzAn{A_n^{\negd*}\negj A_n^\fff}

\newcommand\AnAnz{A_n^\fff\negj A_n^{\negd*}}
\newcommand\AnAnzu{A_n^\fff\negd A_n^{\negd*}}

\newcommand\AkzAk{A_k^{\negd*}\negj A_k^\fff}
\newcommand\BnzBn{B_n^{\negj*}\negj B_n^\fff}

\newcommand\BnBnz{B_n^\fff\negj B_n^*}

\newcommand\AnCnXDnBn{A_n^\fff\negj C_n^\fff\negj X\negj D_{\negj n}^\fff\negj B_{\negj n}^\fff}

\newcommand\Adirei{A_{\dir,\ei}}
\newcommand\Bdirzi{B_{\negj\dir,\zi}}

\newcommand\Astareta{A_{\star,\eta}}

\newcommand\Bstarti{B_{\star,\ti}}

\newcommand\Cdireta{C_{\negj\dir,\eta}}

\newcommand\Dstarti{D_{\!\star,\ti}}

\newcommand\XX{{\mathcal X}}

\newcommand\limm{\lim_{n\ttt\infty}}

\newcommand\jlim[1]{\sideset{^{_{{\sstyle s}\negd}}}\fff\lim_{#1}}
\newcommand\jlimm{\sideset{^{_{{\sstyle s}\negc}}}\fff\limm}

\newcommand\Xlim[1]{\sideset{^{_{\negj{\scriptscriptstyle\XX}\!}}}\fff\lim_{#1}}

\newcommand\wlimm{\sideset{^{_{{\scriptstyle w}\negc}}}\fff\lim_{\! n\ttt\iii}}

\newcommand\nlin\newline
\newcommand\npag\newpage

\newcommand\bskip\bigskip
\newcommand\mmskip\medskip
\newcommand\sskip\smallskip
\newcommand\ima\exists

\newcommand{\mexp}[1]{{\negj\ssstyle(\negd{#1}\negt)}}

\newcommand\eqreff[1]{Eq.\ \eqref{#1}}
\newcommand\ineqreff[1]{Ineq.\ \eqref{#1}}

\newcommand\nji{{n=1}^\infty}

\newcommand\noi{{n=0}^\infty}

\newcommand\uunoi{{\negd n=0}^{\negd\infty}} 

\newcommand\mji{{m=1}^\infty}

\newcommand\mMpji{{m=\MM\plu1}^\iii}

\newcommand\kjn{{k=1}^n}
\newcommand\kon{{k=0}^n}

\newcommand\njN{{n=1}^\NN}

\newcommand\mjM{{m=1}^\MM}

\newcommand\sumN{\sum_{n=1}^\iii}
\newcommand\summ{\sum_{n=1}^\iii}

\newcommand\summud{\summ\negd}
\newcommand\summut{\summ\!}
\newcommand\summuc{\summ\negc}

\newcommand\summm{\sum_{m=1}^\iii}

\newcommand\summmuc{\sum_{m=1}^\iii\negc}

\newcommand\jsummo{\sideset{^{_{\!{\sstyle{s}}\negj}}}\fff\sum_{\negp n=0}^{\negp\iii}}

\newcommand\lwsumm{\sideset{^{_{{\sstyle{w}}}}}\fff{\TS\sum_\nji}}

\newcommand\lwsummo{\sideset{^{_{{\sstyle{w}}\!}}}\fff{\TS\sum_\noi}}

\newcommand\TSBHsummmuc{\sideset{_{\pozx}^{\ssstyle\BH}\negd}\fff{\TS\summm}\negc}

\newcommand\Cssumm{\sideset{^{_{{\ssstyle{\negj\bcalC_{\negd s}\negvi}}}}}\fff\sum_\nji}
\newcommand\Cssummuc{\sideset{^{_{{\ssstyle{\negj\bcalC_{\negd s}\negvi}}}}}\fff\sum_\nji\negc}

\newcommand\Cssummoup{\sideset{^{_{{\ssstyle{\negj\bcalC_{\negd s}\negvi}}}}}\fff\sum_\noi\negp}
\newcommand\TSCssumm{\sideset{^{_{{\ssstyle{\negj\bcalC_{\negd s}\negvi}}}}}\fff{\TS\sum_\nji}}

\newcommand\TSCssummuc{\sideset{^{_{{\ssstyle{\negj\bcalC_{\negd s}\negvi}}}}}\fff{\TS\sum_\nji}\negc}

\newcommand\TSXsumm{\sideset{^{_{{\ssstyle{\negd\XX\negd}}}}}\fff{\TS\sum_\nji}}
\newcommand\TSldlisummmuc{\sideset{^{_{{\ssstyle{\negd\IIS\ldsli\negix}}}}}\fff{\TS\sum_\mji}\negc}
\newcommand\TSldtrosummmuc{\sideset{^{_{{\ssstyle{\negd\IIS\ldstro\negix}}}}}\fff{\TS\sum_\mji}\negc}
\newcommand\TSldlirozsummmuc{\sideset{^{_{{\ssstyle{\negd\IIS\ldsliroz\negix}}}}}\fff{\TS\sum_\mji}\negc}
\newcommand\TSldlirozsummMpjiuc{\sideset{^{_{{\ssstyle{\negd\IIS\ldsliroz\negix}}}}}\fff{\TS\sum_\mMpji}\negc}

\newcommand\jsummut{\negc\sideset{^{_{\sstyle{s}\negc}}}\fff\sum_\nji\!}
\newcommand\jsummuc{\negc\sideset{^{_{\sstyle{s}\negt}}}\fff\sum_\nji\negc}
\newcommand\jsummup{\negc\sideset{^{_{\sstyle{s}\negc}}}\fff\sum_\nji\negp}

\newcommand\TSjsummuc{\sideset{^{_{{\sstyle{s}}\negc}}}\fff{\TS\sum_\nji}\negc}
\newcommand\TSjsummup{\sideset{^{_{{\sstyle{s}}\negc}}}\fff{\TS\sum_\nji}\negp}

\newcommand\TSjsummo{\sideset{^{_{\sstyle{s}}}\negd}\fff{\TS\sum_\noi}}

\renewcommand\XX{{\mathcal X}}                                

\newcommand\HH{{\mathcal H}}

\newcommand\BH{ {\mathcalb B}\negj(\negj{\mathcal H}\negj)}

\newcommand\HHN{\HH^{\negd\NN}}

\newcommand\uuBHHN{{\negd\Be\negd\oduu{\HHN}}}
\newcommand\uuBHHNH{{\negd\Be\negd\oddj{\HHN\negp,\negd\HH}}}
\newcommand\uuBHHHN{{\negd\Be\negd\oduu{\HH\negj,\negd\HHN}}}

\newcommand\ldsBH{{\ell_{\negd s}^2\!\odjjc{\eN,\!\BH}}}
\newcommand\ldsliBH{{\ell_{\! s}^2\!\odjjc{\eN,\negd\lambda,\!\BH}}}
\newcommand\ldsli{{\ell_{\!s\negj,\negj\lambda}^2}}
\newcommand\ldstliBH{{\tilde{\ell}_{\negd s}^2\!\odjjc{\eN,\negd\lambda,\!\BH}}}
\newcommand\ldstli{{\tilde{\ell}_{\negd s\negj,\lambda}^2}}

\newcommand\ldstroBH{{\tilde{\ell}_{\negd s}^2\!\odjjc{\eN,\!\rho,\!\BH}}}
\newcommand\ldstro{{\tilde{\ell}_{\negd s\negj,\rho}^2}}
\newcommand\ldsliroz{{{\ell}_{\negd s\negj,\negj\li\negj,\negj\rho^*}^2}}
\newcommand\ldslirozNBH{\ldsliroz\negc\odjjt{\eN,\!\BH}}

\newcommand\ccj{{\mathcalb C}_{\negj1\!}\negj(\negj{\CMcal H}\negj)}

\newcommand\ccd{{\mathcalb C}_{\negj2\negj}\negj(\negj{\CMcal H}\negj)}

\newcommand\cci{{\mathcalb C}_\iii\negj(\negj{\CMcal H}\negj)}

\newcommand\ccp{{\mathcalb C}_{p\negj}(\negj{\mathcal H}\negj)}

\newcommand\ccs{ {\mathcalb C}_{\negd s\negj}(\negd{\mathcal H}\negj)}

\newcommand\ccphi{ {\mathcalb C}_{\negd\Phi}\negd(\negj{\mathcal H}\negj)}

\newcommand\ccf{ {\mathcalb C}_{\negd\Phi}\negd(\negd{\mathcal H}\negj)}

\newcommand\Phiz{{\Phi^{\negd*}}}

\newcommand\intW{\int_\WWW\negd}

\newcommand\dt{\pozd d\mu(\negj t\negj)}

\newcommand\nizjd[1]{\zzmo{\negj{#1}_n\negd}_{n=1}^\iii}

\newcommand\niznjNu[1]{\zzmo{\negj{#1}_n\negj}_{n=1}^\NN}

\newcommand\nizmexp[1]{\odu{{#1}^\mexp{m}}_\mji}

\newcommand\vniz[1]{\zzvo{{#1}_n}_{n=1}^\iii\!}           

\newcommand\AAnN{{\{A_n^\fff\}_{n=1}^\iii}}     
\newcommand\BBnN{{\{B_n^\fff\}_{n=1}^\iii}}

\newcommand\pozj{\:\!}
\newcommand\pozd{\;\!}
\newcommand\pozt{\,}

\newcommand\pozs{\:\,}

\newcommand\pozx{\;\;}

\newcommand\pozxx{\;\;\;\;}

\newcommand\negj{\;\!\!}               
\newcommand\negd{\:\!\!}

\newcommand\negt{\!}
\newcommand\negc{\;\!\!\!}
\newcommand\negp{\:\!\!\!}

\newcommand\negvi{\!\!}
\newcommand\negs{\;\!\!\!\!}

\newcommand\nego{\:\!\!\!\!}
\newcommand\negix{\!\!\!}
\newcommand\negx{\;\!\!\!\!\!}

\newcommand\negxii{\!\!\!\!}

\newcommand\negxiv{\:\!\!\!\!\!\!}

\newcommand\negxx{\:\!\!\!\!\!\!\!\!}

\newcommand\vspd{\vspace{-2mm}}

\newcommand\vspdp{\vspace{-2.5mm}}

\newcommand\vspc{\vspace{-4mm}}

\newcommand\fff{{\phantom{}}}

\newcommand\llp{\left<}
\newcommand\rrp{\right>}

 \newcommand\zza[1]{\left\vert{#1}\right\vert}   
 \newcommand\zzao[1]{\vert{#1}\vert}

 \newcommand\zzaou[1]{\vert{\negj#1\negj}\vert}

 \newcommand\zzaj[1]{\bigl\vert{#1}\bigr\vert}

 \newcommand\zzajod[1]{\bigl\vert{#1}\negd\bigr\vert}

 \newcommand\zzajtd[1]{\bigl\vert\!   {#1}\negd\bigr\vert}
 
 \newcommand\zzad[1]{\Bigl\vert{#1}\Bigr\vert}            

 \newcommand\zzp[1]{\llp{#1}\rrp}
 \newcommand\zzpo[1]{\langle{#1}\rangle}           
 \newcommand\zzpou[1]{\langle{\negj#1\negj}\rangle}           
 
 \newcommand\zzppu[1]{\zzp{\negc\zzp{#1}\negc}}          
 
 \newcommand\LLnd[1]{\Bigl\vert\!\Bigl\vert{#1}\Bigr.\Bigr.}
 \newcommand\RRnd[1]{\Bigl.\Bigl.{#1}\Bigr\vert\!\Bigr\vert}

 \newcommand\RRndt[1]{\Bigl.\Bigl.{#1}\!\Bigr\vert\!\Bigr\vert}
 \newcommand\LLndc[1]{\Bigl\vert\!\Bigl\vert\negc{#1}\Bigr.\Bigr.}

 \newcommand\zzmo[1]{({#1})}

 \newcommand\zzmou[1]{(\negj{#1}\negj)}

 \newcommand\zzmojd[1]{(\negj{#1}\negd)}

 \newcommand\zzmotd[1]{({#1}\negd)}

\newcommand\odod[1]{(#1\negd)}             
\providecommand\odu[1]{\zzmou{#1}}

\newcommand\odjd[1]{(\negj#1\negd)}
\newcommand\oddj[1]{(\negd#1\negj)}

\newcommand\oduu[1]{(\negd#1\negd)}

\newcommand\oddt[1]{(\negd#1\!)}

 \newcommand\zzmj[1]{\bigl({#1}\bigr)}              

 \newcommand\zzmju[1]{\bigl(\negj{#1}\negj\bigr)}
 
 \newcommand\zzmjod[1]{\bigl(     {#1}\negd\bigr)}
 \newcommand\zzmjjj[1]{\bigl(\negj{#1}\negj\bigr)}

 \newcommand\zzmjdj[1]{\bigl(\negd{#1}\negj\bigr)}

 \newcommand\zzmjuu[1]{\bigl(\negd{#1}\negd\bigr)}

 \newcommand\odju[1]{\bigl(\negj{#1}\negj\bigr)}

 \newcommand\odjjd[1]{\bigl(\negj{#1}\negd\bigr)}

 \newcommand\odjjt[1]{\bigl(\negj{#1}   \!\bigr)}
 \newcommand\odjuu[1]{\bigl(\negd{#1}\negd\bigr)}

 \newcommand\odjte[1]{\bigl(\!   {#1}\pozj\bigr)}
 \newcommand\odjjc[1]{\bigl(\negj{#1}\negc\bigr)}
 \newcommand\odjdt[1]{\bigl(\negd{#1}\!   \bigr)}

 \newcommand\odjut[1]{\bigl(\!   {#1}\!   \bigr)}

\newcommand\uuodjjd[1]{\negd\odjjd{#1}}

 \newcommand\zzmduu[1]{\Bigl(\negd{#1}\negd\Bigr)}
 \newcommand\zzmdtj[1]{\Bigl(\!{#1}\negj\Bigr)}

 \newcommand\zzmdjc[1]{\Bigl(\negj{#1}\negc\Bigr)}
 \newcommand\zzmddt[1]{\Bigl(\negd{#1}\!\Bigr)}
 \newcommand\zzmdtd[1]{\Bigl(\negd{#1}\negd\Bigr)}
 \newcommand\zzmdcj[1]{\Bigl(\negc{#1}\negj\Bigr)}
 
 \newcommand\zzmdjp[1]{\Bigl(\negj{#1}\negp\Bigr)}
 
 \newcommand\zzmddc[1]{\Bigl(\negd{#1}\negc\Bigr)}
 
 \newcommand\zzmdpc[1]{\Bigl(\negp{#1}\negc\Bigr)}

\newcommand\zzmdcd[1]{\Bigl(\negc{#1}\negd\Bigr)}
\newcommand\zzmddp[1]{\Bigl(\negd{#1}\negp\Bigr)}
\newcommand\zzmdtc[1]{\Bigl(\!   {#1}\negc\Bigr)}
\newcommand\zzmdct[1]{\Bigl(\negc{#1}\!   \Bigr)}

\newcommand\zzmduc[1]{\Bigl(\negc{#1}\negc\Bigr)}

\newcommand\zzso[1]{[{#1}]}                  

\newcommand\zzst[1]{\biggl[{#1}\biggr]}

\newcommand\zzv[1]{\left\{{#1}\right\}}         

\newcommand\zzvo[1]{\{{#1}\}}

\newcommand\zzvou[1]{\{\negj{#1}\negj\}}

\newcommand\zzvojd[1]{\{\negj{#1}\negd\}}

\newcommand\zzvj[1]{\bigl\{{#1}\bigr\}}              
\newcommand\zzvju[1]{\bigl\{\negj{#1}\negj\bigr\}}

\newcommand\zzvjdp[1]{\bigl\{\negd{#1}\negp\bigr\}}

\newcommand\zzvjuc[1]{\bigl\{\negc{#1}\negc\bigr\}}

\newcommand\zzvdu[1]{\Bigl\{{\negj#1\negj}\Bigr\}}    

\newcommand\zzvtuu[1]{\biggl\{{\negd#1\negd}\biggr\}}

 \newcommand\nnnou[1]{\zza{\:\!\!\zza{\negj#1\negj}\:\!\!}}

 \newcommand\nnnoo[1]{\zzao{\:\!\!\zzao{#1}\:\!\!}}
 \newcommand\nnnooj[1]{\zzao{\:\!\!\zzao{    #1\negj}\:\!\!}}
 \newcommand\nnnoooj[1]{\zzao{\:\!\!\zzao{    #1\negj}\:\!\!}}
 
 \newcommand\nnnoou[1]{\zzao{\:\!\!\zzao{\negj#1\negj}\:\!\!}}
 
 \newcommand\nnnoojd[1]{\zzao{\:\!\!\zzao{\negj#1\negd}\:\!\!}}

 \newcommand\nnnojod[1]{\zzaj{\!\zzaj{     #1\negj}\!}}
 \newcommand\nnnoju[1]{\zzaj{\!\zzaj{\negj #1\negj}\!}}
 
 \newcommand\nnnojjd[1]{\zzaj{\!\zzaj{\negj#1\negd}\!}}
 
 \newcommand\nnnojjt[1]{\zzaj{\!\zzaj{\negj#1   \!}\!}}
 \newcommand\nnnojuu[1]{\zzaj{\!\zzaj{\negd#1\negd}\!}}
 \newcommand\nnnojtj[1]{\zzaj{\!\zzaj{   \!#1\negj}\!}}
 \newcommand\nnnojtd[1]{\zzaj{\!\zzaj{   \!#1\negd}\!}}
 \newcommand\nnnojdt[1]{\zzaj{\!\zzaj{\negd#1   \!}\!}}
 
 \newcommand\nnnojut[1]{\zzaj{\!\zzaj{   \!#1   \!}\!}}
 \newcommand\nnnojcd[1]{\zzaj{\!\zzaj{\negc#1\negd}\!}}
 \newcommand\nnnojpj[1]{\zzaj{\!\zzaj{\negp#1\negj}\!}}
 
 \newcommand\nnodod[1]{\zzad{\!\zzad{#1\negd}\!}}
 \newcommand\nnodu[1]{\zzad{\!\zzad{\negj#1\negj}\!}}

 \newcommand\nnodjd[1]{\zzad{\!\zzad{\negj#1\negd}\!}}

 \newcommand\nnodjt[1]{\zzad{\!\zzad{\negj#1\!}\!}}
 \newcommand\nnoduu[1]{\zzad{\!\zzad{\negd#1\negd}\!}}
 \newcommand\nnodtj[1]{\zzad{\!\zzad{\!#1\negj}\!}}
 
 \newcommand\nnoddt[1]{\zzad{\!\zzad{\negd#1\!}\!}}
 \newcommand\nnodtd[1]{\zzad{\!\zzad{\!#1\negd}\!}}
 \newcommand\nnodcj[1]{\zzad{\!\zzad{\negc#1\negj}\!}}
 
 \newcommand\nnodut[1]{\zzad{\!\zzad{\!#1\!}\!}}
 \newcommand\nnodcd[1]{\zzad{\!\zzad{\negc#1\negd}\!}}

 \newcommand\nnodct[1]{\zzad{\!\zzad{\negc#1\!}\!}}
 
 \newcommand\nnoduc[1]{\zzad{\!\zzad{\negc#1\negc}\!}}

 \newcommand\nnodsd[1]{\zzad{\!\zzad{\negs#1\negd}\!}}

\newcommand\nnnuo[1]{\zzao{\;\!\!\zzao{\;\!\!\zzao{#1}\;\!\!}\;\!\!}}

\newcommand\matjt[1]{{\zzst{\!\begin{array}{c} #1 \end{array}\!} }}

\newcommand\sbinom\tbinom                                     %
\newcommand\skbinom\tbinom
\newcommand\dsbinom\tbinom

\newcommand\tbinomu[2]{\tbinom{\negj#1\negj}{\negj#2\negj}}

\newcommand\awss{a.w.s.s.\ }
\newcommand\daw{d.a.w.\ }

\newcommand\ONB{orthonormal basis\ }

\newcommand{\Be}{{\mathcalb B}}

\newcommand\N{{\scriptstyle N}}   

\newcommand\NN{{\ssstyle N}}

\newcommand\M{{\scriptstyle M}}
\newcommand\MM{{\ssstyle M}}

\newcommand\Dj{\mbox{\rule[.75ex]{.4em}{.7pt}\kern-.4em D}}
\newcommand\djm{\mbox{\kern.1em\rule[1.2ex]{.3em}{.7pt}\kern-.4em d}}

\newcommand\bfc{{\boldsymbol c}}

\newcommand\bbS{\mathbb{S}}

\newcommand\calA{{\mathcal A}}   

\newcommand\bcalC{{\mathcalb C}}

\newcommand\ai{\alpha}
\newcommand\bi{\beta}
\newcommand\gi{\gamma}

\newcommand\ei{\varepsilon}

\newcommand\zi{\zeta}
\newcommand\ti{\theta}

\newcommand\li{\lambda}

\newcommand\si{\sigma}

\newcommand\WWW{\Omega}

\newcommand\reduk\upharpoonright

\newcommand\Ac{A_{\negj_C}}

\newcommand\emj{{\minjd1}}

\newcommand\ind[1]{{\negj_{\IS#1}}}

\newcommand\utz{{\!*}}
\newcommand\ucz{{\negc*}}

\newcommand\zoxuu{^*\oxuu }
\newcommand\zoxtj{^*\oxtj }
\newcommand\zoxtd{^*\oxtd }

\newcommand\jsummNup{\negc\sideset{^{_{\sstyle{s}\negc}}}\fff\sum_{k=\NN+1}^\iii\negp}

\newcommand\TSjsummNuc{\sideset{^{_{{\sstyle{s}}\negc}}}\fff{\TS\sum_{k=\NN\negj+\negj1}^\iii}\negc}

\newcommand\AkAkz{A_k^\fff\negj A_k^{\negd*}}
\newcommand\BkzBk{B_k^{\negj*}\negj B_k^\fff}
\newcommand\BkBkz{B_k^\fff\negj B_k^*}
\newcommand\Cslimm{\sideset{^{_{{\ssstyle{\negj\bcalC_{\negd s}\negvi}}}}}\fff\lim_{n\to\infty}}

\begin{document}
\def\bibsection{\section*{References}}


\title{\textbf{Extended Cauchy-Schwarz inequalities for $\si$-elementary transformers in Schatten--von
         Neumann norm ideals}}

\author{Danko R.\ Joci\' c$^1$, Mihailo Krsti\' c$^2$, Milan Lazarevi\' c$^3$, Stevan Mila\v sinovi\' c$^4$ \\ \\
$^1$\footnotesize{University of Belgrade, Faculty of Mathematics, Studentski trg 16, Belgrade, Serbia} \\
\footnotesize{\small \texttt{E-mail: danko.jocic@matf.bg.ac.rs}} \\
$^2$\footnotesize{University of Belgrade, Faculty of Mathematics, Studentski trg 16, Belgrade, Serbia}\\
\footnotesize{\small \texttt{E-mail: mihailo.krstic@matf.bg.ac.rs}} \\
$^3$\footnotesize{University of Belgrade, Faculty of Mathematics, Studentski trg 16, Belgrade, Serbia}\\
\footnotesize{\small \texttt{E-mail: milan.lazarevic@matf.bg.ac.rs}}\\
$^4$\footnotesize{University of Belgrade, Faculty of Transport and Traffic Engineering, Vojvode Stepe 305, Belgrade, Serbia}\\
\footnotesize{\small \texttt{E-mail: s.milasinovic@sf.bg.ac.rs}} }

\maketitle

\begin{abstract}
 Let $q,r,s\gesuu1\!$ satisfy $\!\frac1{2q}\plu\frac1{2r}\equt\frac1s\!$ and $X\inuu\ccs.$
If  $\!\nizjd{\lambda},\!\nizjd{w}\!\!$ are sequences in $\!\eRpo$ and $\nizjd{\li_n^{\negc\odu{1\minu q}\negd/\negd\odu{2q}}\negc A},\odjd{\li_n^{\negc1\negd/\negd\odu{2q}}\negc\Anuz}_\nji, \nizjd{w_n^{\mindt1\negd/\negd\odu{2r}}\!B}$ and $\odjd{w_n^{\negd\odu{r\minu1}\negd/\negd\odu{2r}}\!\Bnuz}_\nji\!$ are strongly square summable, then there exists  $\TSCssummuc\AnXBn$ and
\vspbbalzd
\nnoduu{\negd\Cssummuc\AnXBn}_s
\lessd\nnodod{\jsummuc\li_n^{\utfrac1q}\AnAnz}^{\utfrac12\minu\frac1{2q}}\!
\nnodjd{\jsummup w_n^{\mintj\frac1r}\!\BnzBn}^{\utfrac12\minu\frac1{2r}}\!
\nnodct{\zzmddt{\jsummup\li_n^{\utfrac1q\minu1}\negp\AnzAn}^\uvifrac1{2q}\negp X\!
\zzmddc{\jsummup w_n^{\negc1\minu\frac1r}\!\BnBnz}^\uvifrac1{2r}}_s\!.
\end{align*}\vsplintp
Equivalent inequalities are also given, together with some applications to families $\!\nizjd{A}$
and $\!\nizjd{B}$ in $\BH$ which are not double square summable. The results presented in this article
significantly extends the previous results of authors related to $\si$-elementary transformers in Schatten--von Neumann ideals.
\end{abstract}



\tableofcontents
\bigskip

\section{Introduction}

Let $\HH$ be a separable Hilbert space and with $\BH$ we will denote the algebra of all bounded operators on $\HH$. Also, let $\ccphi$ be a ideal of compact operators in $\BH$ generated by s.n.\ function $\Phi$ and the norm of this ideal will be denoted by $\nnnoo{\cdot}_\Phi$. Ideal $\ccphi$ is a Banach space with the norm $\nnnoo{\cdot}_\Phi$. For the theory of normed ideals see \cite{GK} and \cite{Si}. Corresponding ideals of nuclear, Hilbert-Schmidt and compact operators are denoted by $\ccj$, $\ccd$ and $\cci$ with their norms $\nnnoo{\cdot}_1$, $\nnnoo{\cdot}_2$ and $\nnnoo{\cdot}_\iii$ denoted respectively. For the operator norm on the space $\BH$ we will use the symbol $\nnnoo{\cdot}.$ 

The first version of the Cauchy--Schwarz (C--S) norm inequality for $\si$-elementary transformers (of the form) $X\mto\summut\AnXBn$ acting on Schatten--von Neumann ideals $\ccp$ appeared in \cite[Th.~2.1]{J99}.
It is followed by the upgraded version in \cite[Th.~3.3]{J05}, which also includes its generalization to the wider class of inner product type (i.p.t.) transformers $X\mto\intW\AtXBt\dt,$ based on Gel'fand integration of the operator valued functions (operator fields). Some further refinements and alternative and equivalent forms are presented in \cite[Th.~2.1,Th.~2.2]{J23}, contributing to the wider applicability of those inequalities. Various classes of more general transformers, including $\si$-elementary transformers, can be approximated by elementary transformers, see, e.g., \cite[Chapter 5]{AM}.
Therefore it appears well worth to continue to study this class of transformers.

The remarkable aspect of the C--S norm inequalities for $\si$-elementary transformers in Schatten--von Neumann ideals is that they do not require any commutativity properties for the involved families $\AAnN$ and $\BBnN.$ This is in the strong contrast to the C--S norm inequalities for arbitrary unitarily invariant norms (see \cite[Th.~2.2]{J98} and \cite[Th.~3.2]{J05}) and to the C--S norm inequalities in Q and Q$^*$ ideals presented in \cite[Lemma~3.4]{JMDj} and \cite[Th.~3.1]{JKL}, where the normality and commutativity of at least one of those families is required. Such requirement automatically exclude
asymmetrically double square summable sequences from applicability of those C--S norm inequalities.

The main goal of this paper is to provide the extensions of results obtained above in the context of
$\si$--elementary transformers induced by asymmetrically double square summable sequences $\nizjd{A}$
and $\nizjd{B}\negj;$ see Notation subsection in this paper for the precise definitions.

Recall that for selfadjoint $A,B\inu\BH$ it is said that $A\lesu B$ if and only if $\zzpo{Ah,h}\lesu\zzpo{Bh,h}$ for all $h\inu\HH.$
Also, if $J\sset\eR$ is an interval and $f$ be a real-valued function on $J,$
it is said that $f$ is operator monotone if, for every selfadjoint $A,B\inu\BH$ with $\si\odu{A}\cupu\si\odu{B}\sset J,$ the operator inequality $A\lesu B$ implies
$f\odu{A}\lesuu f\odu{B},$ where $f\odu{A},f\odu{B}$ are understood in the sense of spectral theorem calculus. For $p\inu\OJ$ the functions $t\mto t^p$ on $\eRp$ are the best known examples of operator monotone functions on $\eRp;$ see \cite[Def.~2.1.3(1), Ex.~2.5.9(3)]{Hi}.

If $A,B\inu\BH$ and $0\lesu A\lesu B$ then we have monotonicity of operator norm $\nnnoo{\cdot}$ i.e.\ $\nnnoo{A}\lesu\nnnoo{B}.$

In our work we will use the following double monotonicity property \eqref{duplamonotonost}, saying that
\vspbbald\label{duplamonotonost}
\nnnoo{AXB}_\ind\Phi\lestj\nnnoo{CXD}_\ind\Phi\quad\text{for all $X\inu\ccphi,$}
\end{align}\vsplincp
whenever $A^*A\lesdj C^*C$ and $BB^*\lesdj DD^*$ (for the proof see \cite[p.~62]{JLM}).

Each $\nnnoo{\cdot}_\Phi$ is lower semicontinuous, i.e., for every weakly convergent sequence $\vniz{X}$ in $\ccf$
\vspbbald\label{poluneprekidnostodozdo}
\nnnojpj{\wlimm X_n}_\ind\Phi\lestj\liminf_{n\to\iii}\nnnoooj{X_n}_\ind\Phi.
\end{align}\vsplinp
This follows from the well-known formula
\vspbbalzd
\nnnoo{X}_\ind\Phi\equu\sup\zzvtuu{\frac{|\tr(XY)|}{\nnnoo{Y}_\ind{\Phi^*}}:
\text{$Y$ is a finite rank nonzero operator}},
\end{align*}\vsplinc
where the finite rank operators are those of the form
$\sum_{k=1}^ng_k\zoxuu h_k\colon\HH\ttt\HH\colon f\mto\sum_{k=1}^n\zzpo{f,g_k}h_k$ and
$\nnnoo{\cdot}_\ind\Phiz\!\!$ denote the dual norm for $\nnnoo{\cdot}_\ind\Phi.$
Here, for any $f,g\inu\HH,$ we denote by $g\zoxuu f$ one dimensional operator acting by
$g\zoxuu f(h)\eqdu\zzpo{h,g} f$ for any $h\inu\HH,$ which are known to have their
linear span dense in each of $\ccp$ for any $1\lesu p\lesu\iip.$

\subsection{Notations and terminology}
In the sequel $\eR$ will stand for real numbers, $\Rp\eqdtj\eRp,\Rpz\eqdtj\eRpz$ and $\eN$ will denote
the set of positive integers.

$\bbSOJoe$ (resp.\ $\bbS\OJ$) will denote the open vertical strip $\bbSOJo\eqdu\zzvou{z\inuu\Ce:0<\Re z<1}$ (resp.\ closed vertical strip $\bbS\OJ\eqdu\zzvou{z\inuu\Ce:0\lesuu\Re z\lesuu1}$).

In a normed space  $(\XX,\nnnoo{\cdot})$  and a sequence $\vniz{\boldsymbol{x}}$ we  write
$\boldsymbol{x}=\Xlim{\negd n\to\iii}{\boldsymbol{x}}_n$ for some $\boldsymbol{x}\in\XX$ if $\lim_\nti\nnnoo{{\boldsymbol{x}}_n-\boldsymbol{x}}=0.$ Also, $\boldsymbol{x}=\TSXsumm\boldsymbol{x}_n$
will denote that $\boldsymbol{x}=\Xlim{\negd n\to\iii}\sum_\kjn\boldsymbol{x}_k.$ If $\XX\dod\ccs$ for
the Schatten $\ccs$ ideal for some $1\lesu s\lesu\iip,$ we use the abbreviation $\!\Cslimm$ (resp.\
$\negp\TSCssumm$) for the limits (resp.\ sums) in the Schatten $\ccs$ ideals.
The exceptions from this notation will be the strong (resp.\ weak) operator
limit and sum $\!\jlimm\!$ and $\!\TSjsummo$ (resp.\ $\negc\wlimm\!$ and $\!\lwsummo$\!).

We refer to a sequence $\nizjd{A}$ in $\BH$  as a (strongly) \textit{square summable} (s.s.) if 
$\summuc\nnnoo{A_n h}^2<\iip$ for all $h\inu\HH.$
By using the polarization identity this imply that $\zzv{\sum_\kjn\!\AkzAk}_\nji$ weakly
converges to $\lwsumm\negp\AnzAn\inu\BH$ as $n\ttt\iip,$ according to the uniform boundedness
principle. Due to the monotonicity of its partial sums $\sum_\kjn\!\AkzAk,$ the considered
convergence is moreover strong. If there exists both $\TSjsummuc\AnzAn$ and $\TSjsummuc\AnAnz,$ then
we say that $\nizjd{A}$ is \textit{double square summable} (d.s.s).
Moreover, if there exists $\TSjsummuc w_n^\fff\AnzAn$ for some (weight)
sequence $\nizjd{w}$ in $\Rpz,$ then $\nizjd{A}$ will be called \textit{$w$--square summable} (w-s.s.) sequence in $\BH.$

For $\N$-tuples $\odu{A_n}_\njN$ in $\BH$ we use notations $[A_n]_\njN$ for the $1\x\N$ row operator matrix $[A_1,\ldots,A_{\negd\NN}]\inu\Be\odju{\HH^{\negd\NN}\negc,\!\HH}.$
For typographical reasons it is easier to write row vectors, so we will write
$\odju{[A_n]_\njN}^{\!\!\top}\eqdu [A_1,\ldots,A_{\negd\NN}]^{\!\top}$ for a typical $\N\x1$ column operator matrix
$\matjt{\nego\IS A_1\vspdp\negxii\\\nego\IS \vdots\vspd\negxii\\ \nego\IS A_{\negd\NN}\negxii\vspace{1mm}}\inu\Be\uuodjjd{\HH,\!\HH^{\negd\NN}},$
where, $\IS"\top"$ denotes the transpose.
Let us note that according to $C^*\!$ property for operators acting between Hilbert spaces
\vspbbald\label{Czvezda}
&\nnnojuu{[\Anff]_\njN}^2_\uuBHHNH\eqcj\nnnojtd{[\Anff]_\njN\!\odjuu{[\Anff]_\njN}^\ucz}
=\nnnojut{[\Anff]_\njN\!\odjuu{[A_n^{\negd*}]_\njN}^{\!\!\top}}=\nnoduu{\sum_\njN\AnAnz},\\
&\nnnojut{\odjuu{[\Anff]_\njN}^{\!\!\top}}_\uuBHHHN^2
\eqcj\nnnojut{\odjuu{\odjuu{[\Anff]_\njN}^{\!\!\top}}^\ucz\negc\odjuu{[\Anff]_\njN}^{\!\!\top}}
=\nnnojut{[\Anuz]_\njN\!\odjuu{[A_n]_\njN}^{\!\!\top}}=\nnoduu{\sum_\njN\AnzAn}.\label{CzvezdaAlt}
\end{align}\vsplincp

For $C,D\inu\BH$ the bilateral multiplier $C\oxdt D$ is defined as the transformer
$C\oxdt D:\BH\to\BH:X\mto CXD,$ and $\GCD$ stands for $I\oxdt I\minu C\oxdt D:\BH\to\BH:X\mto X\minuu CXD.$

If $\N\inN,$ then $C\inu\BH$ is \textit{$\N$-hypercontraction} if
$\GCzCI{n}\eqdc\sum_\kon\negj\odu{\minu1}^k\negd\tbinom{n}k C^{*k} C^k\gesu0$ for all $1\lesu n\lesu\N,$
and similarly, $C\inu\BH$ is \textit{$\N$-cohypercontraction} if
$\GCCzI{n}\eqdc\sum_\kon\negj\odu{\minu1}^k\negd\tbinom{n}k C^k C^{*k}\gesu0$ for all $1\lesu n\lesu\N.$

We use the notation $\Ac\eqdu\calA(C)\eqdu\jlimm C^{*n}C^n\negd$ for every contractive $C\inu\BH.$
If $\Ac=0,$ then a contraction $C$ and its semigroup $\zzvojd{C^n}_\noi$ are called \textit{stable.}
For the existence, additional properties and the role of $\Ac$ for hypercontractive operators and in the stability of operators and operator semigroups, see \cite[Prop.~3.1]{K}, \cite[Lemma~2.7]{JL},
\cite[Th.~2.23, Cor.~2.25]{E}, and the references therein.

We also need to emphasize that throughout this paper we will treat (address to) every unnumbered line
in a multiline formula as (to) a part of the consequent numbered one.

\section{Preliminaries}
\subsection{Double modules of asymmetrically weighted square summable operator families}

\bbdef
If $\li\dodu\nizjd{\li}$ is a sequence in $\Rpz,$ let
\bbn
\item $\ldsliBH$ (and $\ldsli$ for short) to denote the class of all $\li$-s.s.\ sequences $\nizjd{A}$ in $\BH,$ i.e.\ those that there exists $\TSjsummuc\li_n\AnzAn\odjjt{\inu\BH},$
\item $\ldstliBH$ (and $\ldstli$ for short) to denote the class of all sequences $\nizjd{A}$ in $\BH$ such that there exists $\TSjsummuc\li_n\AnAnz,$
\item if $\li\dodu\mathbf{1},$  i.e.\ if $\li_n\dodu1$ for all $n\inN,$
then we will also use the simplified notations $\ldsBH$ and $\ell^2\!\odjut{\BH}$ for $\ldsliBH.$
\een
\eedef
It is well known that $\ldsliBH$ (including $\ldsBH$) is the right Hilbert $C^*\!$-module with the inner product and the norm given by
\vspbbald
\zzppu{A,B}_\ind{\!\li}&\eqdto\jsummuc\li_n^\fff\Anuz B_n^\fff\qquad \textrm{for any $A,B\inu\ldsliBH,$}\\
\nnnou{A}_\ind{\negj\li}&\eqdto\nnodjd{\jsummuc\li_n^\fff\AnzAn}^\uujd\quad \textrm{for any $A\inu\ldsliBH.$}
\end{align}\vsplinc
For more details see also \cite[p.p.~6]{L} and \cite[Ex.~2.5.5]{MT}. 
Similarly, $\ldstliBH$ is the left $C^*\!$-module with the inner product and the norm given by
\vspbbald
\zzppu{A,B}_\ind{\negc\li^{\negd*}}&\eqdto\jsummuc\li_n^\fff A_n^\fff B_n^*\qquad \textrm{for any $A,B\inu\ldstliBH,$}\\
\nnnou{A}_\ind{\negj\li^{\negd*}}&\eqdto\nnodjd{\jsummuc\li_n^\fff\AnAnz}^\uujd\quad \textrm{for any $A\inu\ldstliBH.$}
\end{align}\vsplins
\bbdef
$\ldsliroz\negc\eqdu\ldsliroz\negc\odjjt{\eN,\!\BH}\eqdu\ldsliBH\captj\ldstroBH$ will be referred
as a \textit{double, asymmetrically weighted} (d.a.w.) (Hilbert) C$^*\!$-module.

Any operator sequence $\nizjd{A}\inu\ldsliroz$ will be referred as a \textit{$\odod{\li,\rho^*}$\!\! -asymmetrically weighted square summable} $\odjut{\odod{\li,\rho^*}\minct\textit{a.w.s.s.}},$ while for $\odod{\mathbf{1},\rho^*}$-a.w.s.s.\ sequence of operators we simply say that it is a
$\rho^*\negj$-\awss sequence.
\eedef
These spaces are considered in the following lemma.
\bbl
If $\li\dodu\nizjd{\li}$ and $\rho\dodu\nizjd{\rho}$ are sequences in $\Rpz,$ then \daw Hilbert
C$^*\!$-module $\ldslirozNBH$ is a complete if its  norm is introduced by
\vspbbald\label{MaxNorm}
\nnnou{A}_\ind{\negj\li,\rho^*}&\eqdu\max\zzvdu{\nnodjd{\jsummuc\li_n^\fff\AnzAn}^\uujd\!,\!
\nnodjd{\jsummuc\rho_n^\fff\AnAnz}^\uujd}\quad \textrm{for any $A\inu\ldslirozNBH.$}
\end{align}\vsplinx
\eel
\bbproo
 To prove this lemma, we start with an arbitrary absolute convergent sequence $\nizmexp{A}$ in $\ldslirozNBH,$ i.e.\ satisfying $\bfc\dod\summmuc\nnnoou{A^\mexp{m}}_\ind{\negj\li,\rho^*}\lsuu\iip.$ Since $\kor{\li_n}\summmuc\nnnou{A_n^\mexp{m}}\lesu\summmuc\nnnou{A^\mexp{m}}_\ind{\negj\li}
\lespt\summmuc\nnnoou{A^\mexp{m}}_\ind{\negj\li,\rho^*}\equ\bfc$ for any $n\inu\eN,$ it follows that there exist $S_n\dod\!\!\TSBHsummmuc A_n^\mexp{m}$ for any $n\inu\eN,$ as $\BH$ is a Banach space, and there is also $S\dodu\TSldlisummmuc A^\mexp{m}\!,$ as $\ldsliBH$ is a Hilbert C$^*\!$-module. Note also that $S=\nizjd{S}$ and $\nnnoo{S}_\ind{\negj\li}=\nnnoju{\TSldlisummmuc A^\mexp{m}}_\ind{\negj\li}
\lesuu\bfc.$
Similarly, $\kor{\rho_n}\summmuc\nnnou{A_n^{\mexp{m}*}}\lesu\summmuc\nnnou{A^\mexp{m}}_\ind{\negd\rho^*}
\lespt\summmuc\nnnoou{A^\mexp{m}}_\ind{\negj\li,\rho^*}\equ\bfc$ for any $n\inu\eN$ implies that there exist $T_n\dod\!\!\TSBHsummmuc A_n^{\mexp{m}*}=\odjut{\TSBHsummmuc A_n^\mexp{m}}^{\!*}\eqcj S_n^*$ for any $n\inu\eN,$ and moreover, there exists $\tilde{S}\dodu\!\TSldtrosummmuc A^{\mexp{m}}\!,$ also satisfying $\nnnooj{\tilde{S}}_\ind{\negd\rho^*}\eqxd\nnnojuu{\TSldtrosummmuc A^\mexp{m}}_\ind{\negd\rho^*}\lesxd \bfc.$ Furthermore, $\tilde{S}=\odu{T_n^*}_\nji\eqdj\nizjd{S}\eqdj S,$ which shows $S\inu\ldstroBH$ as well, and so 
$S\inu\ldslirozNBH$ $=\ldsliBH\captj\ldstroBH,$ additionally satisfying the estimate
\vspbbaljp
\nnnoo{S}\ind{\negj\li,\rho^*}=\max\zzvjdp{\nnnoo{S}_\ind{\negj\li},\nnnoo{S}_\ind{\negj\rho^*}}\lespj\bfc
\eqjt\TS\summmuc\nnnoou{A^\mexp{m}}_\ind{\negj\li,\rho^*}.\label{slabaONzaPresek}
\end{align}\vsplincp
Finally, by applying \ineqreff{slabaONzaPresek} to 
$\!\TSldlirozsummMpjiuc A^\mexp{m}$ instead of $S$ we get
\vspbbald
\nnnojjt{S\minuu\!\TS\sum_\mjM\!A^\mexp{m}}_\ind{\negj\li,\rho^*}
\eqxd\nnnojcd{\TSldlirozsummMpjiuc A^\mexp{m}}_\ind{\negj\li,\rho^*}
\lesxd\TS\sum_{m=\MM\plu1}^\iii\!\nnnojjd{A^\mexp{m}}_\ind{\negj\li,\rho^*}\!\!\!\ttt0
\qquad\textrm{as $\M\ttt\iip,$}\label{slabaONjeJaka}
\end{align}\vsplinc which validates that $S$ is exactly the existing $\TSldlirozsummmuc A^\mexp{m}\!,$ and so $\ldslirozNBH$ is complete.
\eeproo

\bbr
If $\li\dodu\nizjd{\li}\inu\ell^1_{\IIS\eN}$ and $\nnnou{\li}_\ISj\lesu1,$ then $\ldsBH\captj\ldstliBH=\ldsBH,$
as $\nnnou{A_n}^2\lestj\nnnojtd{\TSjsummuc\AnzAn}$ for any $n\inu\eN,$ and therefore
$\nnnojtd{\TSjsummuc\li_n\AnAnz}\lestj\sup_{n\in\eN}\sum_\kjn\zzaou{\li_k}\nnnojjd{\AkAkz}\lestj\nnnou{\li}_\ISj\!\nnnojtd{\TSjsummuc\AnzAn}.$
This shows the existence of the isometric imbedding of $\ldsBH\captj\ldstliBH$ into $\ldsBH.$
\eer
\bbex\label{osnovniprimer}
If $\li\dodu\nizjd{\li}$ in $\Rpz$ is summable and $\nizjd{e}$ is an \ONB in $\HH,$ then
$\odjjd{e_n\zoxtd e_1}_\nji\inu\ldsBH\captj\ldstliBH.$ Indeed,
\vspbbaldp\label{ektremPrimer}
\jsummut\zzajod{e_n\zoxtd e_1}^2\eqdp\jsummut\odjjd{e_1\zoxtd e_n}\negd\odjjd{e_n\zoxtd  e_1}
\eqdp\jsummut e_n\zoxtd e_n\equ I,\\
\jsummut\li_n\zzajtd{\odjjd{e_n\zoxtd  e_1}^\utz}^2
\eqdp\jsummut\li_n\odjjd{e_n\zoxtd e_1}\negd\odjjd{e_1\zoxtd e_n}
\eqdt\zzmdjc{\sumN\li_n} e_1\zoxtd  e_1.
\label{ektremPrimer2}
\end{align}\vsplinxi
\eeex
\bbex\label{dvaesdva}
If $C$ (resp.\ $D$) is $\N$-cohypercontractive (resp.\ $\M$-cohypercontractive) for some $\N,\M\inu\eN,$ then the operator sequence $C_{\!\NN}\dodu\odjte{\TSsqrt{\GCzCI{}\!}\,C^n\!\TSsqrt{\GCCzI{\NN}\!}}_\uunoi$ is $\li^{\negd*}\negd$-\awss for $\li\dodu\zzvjuc{\tbinomu{n\plu\NN\minuu1}{\NN\minuu1}}_\uunoi,$
$D_{\!\MM}\dodu\odjte{\TSsqrt{\GDzDI{}\!}\,D^n\!\TSsqrt{\GDDzI{\MM}\!}}_\uunoi$ is $\rho^*\negd$-\awss for $\rho\dodu\zzvjuc{\tbinomu{n\plu\MM\minuu1}{\MM\minuu1}}_\uunoi,$
i.e.\ $C_{\!\NN}\inu\ldsBH\captj\ldstliBH$ and $D_{\!\MM}\inu\ldsBH\captj\ldstroBH,$ satisfying
\vspbbald\label{hiperkontr1}
\jsummo\!\!\TSsqrt{\GCCzI{\NN}\!}\,C^{*n}\GCzCI{}C^n\!\TSsqrt{\GCCzI{\NN}\!}\,&
=\TSsqrt{\GCCzI{\NN}\!}\,\odjjd{I\minuu\calAC}\negd\TSsqrt{\GCCzI{\NN}\!}\lesod\GCCzI{\NN},\\
\jsummo\!\!\tbinomu{n\plu\NN\minuu1}{\NN\minu1}\negd
\TSsqrt{\GCzCI{}\!}\,C^n\GCCzI{\NN}C^{*n}\!\TSsqrt{\GCzCI{}\!}\,
&=\TSsqrt{\GCzCI{}\!}\,\odjjt{I\minuu\calACz}\negj\TSsqrt{\GCzCI{}\!}\lesod\GCzCI{},\label{hiperkontr2}\\
\jsummo\!\!\TSsqrt{\GDDzI{\MM}\!}\,D^{\negj*n}\GDzDI{}D^n\!\TSsqrt{\GDDzI{\MM}\!}\,&
=\TSsqrt{\GDDzI{\MM}\!}\,\odjjd{I\minuu\calAD}\negd\TSsqrt{\GDDzI{\MM}\!}\lesod\GDDzI{\MM},\label{hiperkontr3}\\
\jsummo\!\!\tbinomu{n\plu\MM\minuu1}{\MM\minuu1}\negd
\TSsqrt{\GDzDI{}\!}\,D^n\GDDzI{\MM}D^{\negj*n}\!\TSsqrt{\GDzDI{}\!}\,
&=\TSsqrt{\GDzDI{}\!}\,\odjjt{I\minuu\calADz}\negj\TSsqrt{\GDzDI{}\!}\lesod\GDzDI{}.\label{hiperkontr4}
\end{align}\vsplinc
Here inequalities \eqref{hiperkontr1} and \eqref{hiperkontr2} follows from Eq.\ (2.10) in \cite[lemma~2.7]{JL} and due to fact that $\calAC$ and $\calACz\negp$ are positive contractive operators,
while \ineqreff{hiperkontr3} and \ineqreff{hiperkontr4} are their duplicates.
\eeex$ $\vsplinx
\bbex
If $C$ is $\N$-hypercontractive and $\M$-cohypercontractive for some $\N,\M\inu\eN,$ then
the operator sequence $C_{\!\MM\negd,\NN}\dodu\odjte{\TSsqrt{\GCzCI{\NN}\!}\,C^n\!\TSsqrt{\GCCzI{\MM}\!}}_\uunoi$ is $\odod{\li,\rho^*}$-\awss for $\li\dodu\zzvjuc{\tbinomu{n\plu\NN\minuu1}{\NN\minuu1}}_\uunoi$ and
$\rho\dodu\zzvjuc{\tbinomu{n\plu\MM\minuu1}{\MM\minuu1}}_\uunoi$,
i.e.\ $C_{\!\MM\negd,\NN}\intj\ldsliBH\cap\ldstroBH,$ as
\vspbbald
\jsummo\!\!\tbinomu{n\plu\NN\minu1}{\NN\minu1}\negd
\TSsqrt{\GCCzI{\MM}\!}\,C^{*n}\GCzCI{\NN}C^n\!\TSsqrt{\GCCzI{\MM}}\lesu\GCCzI{\MM},\label{hiperkontr5}\\
\jsummo\!\!\tbinomu{n\plu\MM\minu1}{\MM\minu1}\negd
\TSsqrt{\GCzCI{\NN}\!}\,C^n\GCCzI{\MM}C^{*n}\!\TSsqrt{\GCzCI{\NN}}\lesu\GCzCI{\NN}.\label{hiperkontr6}
\end{align}\vsplinc
Once again, inequalities \eqref{hiperkontr5} and \eqref{hiperkontr6} follows from Eq.\ (2.10) in \cite[lemma~2.7]{JL}.
\eeex

\subsection{Enhancing Cauchy--Schwarz inequalities for $\si$-elementary transformers in $\BH,\ccj$ and $\ccd$}

In the sequel we will use the following notation.
\bbdef\label{defInvertibilni}
For sequences $\niznjNu{A},\niznjNu{B},\niznjNu{C},\niznjNu{D}\!$ in $\BH,$ $\niznjNu{\li},\niznjNu{w}\!$
in $\Rp\!$ and $\ei,\zi,\eta,\ti$ in $\Rpz$ let
$\Astareta\eqdu\zzmju{\eta I\plujt\sum_\njN\!\li_n\AnAnz}^\uujd\!,
\Bdirzi\eqdu\zzmju{\zi I\pludc\sum_\njN\!\BnzBn}^\uujd\!,$
$\Adirei\eqdu\zzmju{\ei I\pludc\sum_\njN\!\AnzAn}^\uujd\!,$
$\Bstarti\eqdu\zzmju{\ti I\pludc\sum_\njN\!w_n\BnBnz}^\uujd\!,$
$\Cdireta\eqdu\zzmju{\eta I\pludc\sum_\njN\!\CnzCn}^\uujd\!$ and
$\Dstarti\eqdu\zzmju{\ti I\pludc\sum_\njN\!w_n\DnDnz}^\uujd\!.$
\eedef
\bbt\label{extendedCSinCp}
Under conditions of Definition \ref{defInvertibilni} we have
\vspbbald\label{nulaABCD}\,\,
\nnodu{\Astareta^\mj\negc\sum_\njN\!\li_n^\ujd\AnCnXDnBn\Bdirzi^\emj}\negc
&\lesoj\max_{1\lesu n\lesoj\NN}\negc\nnnoou{\Cn}\nnnoou{\Dn}\nnnoou{X}\quad\text{for all}\,\,\, X\inu\BH,\\
\nnodu{\sum_\njN\!A_n^\fff\Cn\Cdireta^\emj X\Dn\Bn\Bdirzi^\emj}_2
&\lespj\max_{1\lesu n\lesoj\NN}\negc\nnnoou{A_n^\fff}\max_{1\lesu n\lesoj\NN}\negc\nnnoou{\Dn}\nnnoou{X}_\ISd\quad\text{for all}\,\,\, X\inu\ccd,\label{dvaABCD2}\\
\nnodu{\Astareta^\mj\!\sum_\njN\!\li_n^\ujd w_n^\ejd A_n^\fff\Cn X\Dstarti^\emj\Dn\Bn}_2
&\lespj\max_{1\lesu n\lesoj\NN}\negc\nnnoou{\Bn}\max_{1\lesu n\lesoj\NN}\negc \nnnoou{\Cn}\nnnoou{X}_\ISd\label{dvaABCD}\quad\text{for all}\,\,\, X\inu\ccd,\\
\nnodu{\sum_\njN\!w_n^\ejd A_n^\fff\Cn\Cdireta^\emj X\Dstarti^\emj\Dn\Bn}_1
&\lespj\max_{1\lesu n\lesoj\NN}\negc \nnnoou{A_n^\fff}\nnnoou{\Bn}\nnnoou{X}_\ISj\label{jenABCD}\quad\text{for all}\,\,\, X\inu\ccj,
\end{align}\vsplinc
which is equivalent to
\vspbbald\label{nulaABCDalt}
\nnodu{\sum_\njN\!\li_n^\ujd\AnCnXDnBn}
&\lesu\max_{1\lesu n\lesoj\NN}\negc\nnnoou{\Cn}\nnnoou{\Dn}
\nnodu{\sum_\njN\!\li_n\AnAnz}^\uujd\negd\nnodu{\sum_\njN\!\BnzBn}^\uujd\negd\nnnoou{X},\\
\nnodu{\sum_\njN\!A_n^\fff\Cn X\Dn\Bn}_2
&\lespj\max_{1\lesu n\lesoj\NN}\negc\nnnoou{A_n\fff}\max_{1\lesu n\lesoj\NN}\negc\nnnoou{\Dn}
\nnodu{\sum_\njN\!\BnzBn}^\uujd\nnodtd{\zzmdjc{\sum_\njN\!\CnzCn}^\upjd\negp X}_\ISd,\label{dvaABCD2alt}\\
\nnodjt{\sum_\njN\!\li_n^\ujd w_n^\ejd\AnCnXDnBn}_2
&\lespj\max_{1\lesu n\lesoj\NN}\negc\nnnoou{\Bn}\max_{1\lesu n\lesoj\NN}\negc\nnnoou{\Cn} \nnodu{\sum_\njN\!\li_n\AnAnz}^\uujd\negd\nnodjd{X\!\zzmdjc{\sum_\njN\!w_n\DnDnz}^\upjd}_\ISd,\label{dvaABCDalt}\\
\nnodu{\sum_\njN\!w_n^\ejd\AnCnXDnBn}_1
&\lespj\max_{1\lesu n\lesoj\NN}\negc\nnnoou{A_n^\fff}\nnnoou{\Bn}
\nnodtd{\zzmdjc{\sum_\njN\!\CnzCn}^\upjd\negp X\!\zzmdjc{\sum_\njN\!w_n\DnDnz}^\upjd}_\ISj\label{jenABCDalt}.
\end{align}\vsplinc
Both inequality quadruple \eqref{nulaABCD}--\eqref{jenABCD} or \eqref{nulaABCDalt}--\eqref{jenABCDalt} implies
\vspbbald
\,\,\nnodu{\Astareta^\mj\!\sum_\njN\!\li_n^\ujd\AnXBn\Bdirzi^\emj}&\lesoj\nnnoou{X}\quad
\text{for all}\,\,\,X\inu\BH,\label{nula}\\
\nnodu{\sum_\njN\!A_n^\fff\Adirei^\mj X\Bn\Bdirzi^\emj}_2&\lespj\nnnoou{X}_\ISd\label{dva2}\quad
\text{for all}\,\,\, X\inu\ccd,\\
\nnodu{\Astareta^\mj\!\sum_\njN\!\li_n^\ujd w_n^\ejd A_n^\fff X\Bstarti^\mj\Bn}_2&\lespj\nnnoou{X}_\ISd\label{dva}\quad\text{for all}\,\,\, X\inu\ccd,\\
\nnodu{\sum_\njN\!w_n^\ejd A_n^\fff\Adirei^\mj X\Bstarti^\emj\Bn}_1&\lespj\nnnoou{X}_\ISj\quad\text{for all}\,\,\, X\inu\ccj\label{jen}.
\end{align}\vsplinix
\eet
\bbproo
Based on the notation from Definition \ref{defInvertibilni}, the proof for \ineqreff{nulaABCD} follows from the estimates
\vspbbald\notag
&\nnodu{\Astareta^\mj\negc\sum_\njN\!\li_n^\uujd\AnCnXDnBn\Bdirzi^\emj}
\equ\nnoduu{\zzso{\li_n^\uujd\negd\Astareta^\mj A_n^\fff}_\njN\!
\bigoplus_\njN\negj\CnXDn\uuodjjd{\zzso{\Bn\Bdirzi^\emj}_\njN}^{\!\!\top}}\\
&\lesoj\nnnojjt{\zzso{\li_n^\uujd\negd\Astareta^\mj A_n^\fff}_\njN}_\uuBHHNH \! \nnodod{\bigoplus_\njN\negj\CnXDn}_\uuBHHN\! \nnnojut{\odjjd{\zzso{\Bn\Bdirzi^\emj}_\njN}^{\!\!\top}}_\uuBHHHN\notag\\
&\eqoj\nnodu{\Astareta^\mj\!\sum_\njN\!\li_n\AnAnz\Astareta^\mj}^\uujd\!\!
\max_{1\lesu n\lesoj\NN}\!\!\nnnoou{\Cn X\Dn}\nnodu{\Bdirzi^\emj\!\sum_\njN\!\BnzBn\Bdirzi^\emj}^\uujd
\lespj\max_{1\lesu n\lesoj\NN}\negc\nnnoou{\Cn}\nnnoou{\Dn}\nnnoou{X},\label{nulaABCDdok}
\end{align}\vsplinc
where the last equality in \eqref{nulaABCDdok} follows by calculating the row and the column norms
according to formulas \eqref{Czvezda} and \eqref{CzvezdaAlt},
while the last inequality in \eqref{nulaABCDdok} is based on the inequalities
$\Astareta^\mj\!\sum_\njN\!\li_n\AnAnzu\Astareta^\mj\lesuu\Astareta^\mj\negj\Astareta^2\Astareta^\mj\equu I$ and (similarly derived) $\Bdirzi^\emj\!\sum_\njN\!\BnzBn\Bdirzi^\emj\lesu I,$ also combined with the monotonicity of the operator norm $\nnnoo{\cdot}$ for positive operators.\smallskip

The proof of \ineqreff{dvaABCD2} is based on
\vspbbald\notag
&\nnodod{\sum_\njN\! A_nC_n\Cdireta^\emj XD_n\Bn\Bdirzi^\emj}_2\\
&\lesod\nnodcd{\zzmdct{\Cdireta^\emj\!\sum_\njN\! C_n^*\AnzAn\Cn\Cdireta^\emj}^\upjd\negp X}_2\negd
\nnodjd{\Bdirzi^\emj\!\sum_\njN\! B_n^*\DnzDn\Bn\Bdirzi^\emj}^\uujd\label{dvaABCD2dok1}\\
&\lespd \nnodcd{\zzmdct{\Cdireta^\emj\!\sum_\njN\ C_n^*\nnnoojd{A_n^\fff}^2\Cn\Cdireta^\emj}^\upjd\negp X}_2\negd
\nnodjd{\Bdirzi^\emj\!\sum_\njN\!B_n^*\nnnoojd{D_n}^2\Bn\Bdirzi^\emj}^\uujd
\lespd\max_{1\lesu n\lesoj\NN}\negc\nnnoojd{A_n^\fff}
\max_{1\lesu n\lesoj\NN}\negc\nnnoojd{\Dn}\nnnoou{X}_\ISd,\label{dvaABCD2dok2}
\end{align}\vsplinc
where \ineqreff{dvaABCD2dok1} is due to the Cauchy-Schwarz norm inequality (10) in \cite[Lemma~2.1(d2)]{JKL}, applied to the families $\zzvj{A_nC_n\Cdireta^\emj}_\njN$ and $\zzvj{D_nB_n\Bdirzi^\emj}_\njN$ instead of
$\AAnN$ and $\BBnN$ therein. The first inequality in \eqref{dvaABCD2dok2} follows from the estimates
$\AnzAn\lesu\nnnoojd{A_n^\fff}^2I$ and $\DnzDn\lesu\nnnoojd{\Dn}^2I,$
while the second inequality in \eqref{dvaABCD2dok2} is based on the operator inequalities
$\Cdireta^\emj\!\sum_\njN\!\CnzCn\Cdireta^\emj\les I$ and $\Bdirzi^\emj\!\sum_\njN\!\BnzBn\Bdirzi^\emj
\les I, $ all combined with the double monotonicity property \eqref{duplamonotonost}. \smallskip

To prove \ineqreff{dvaABCD} we use
\vspbbald\notag
&\nnodjd{\Astareta^\mj\negc\sum_\njN\!\li_n^\ujd w_n^\ejd\! A_nC_nX\Dstarti^\emj D_nB_n}_2\\
&\lesod\nnoduu{\Astareta^\mj\negc\sum_\njN\!\li_n A_nC_nC_n^*A_n^*\Astareta^\mj}^\uujd\negd \nnoduu{X\!\zzmdct{\Dstarti^\emj\negc\sum_\njN\!w_n D_n\BnBnz D_n^*\Dstarti^\emj}^\upjd}_2 \label{Cetvorka1}\\
&\lespd\nnoduu{\Astareta^\mj\negc\sum_\njN\!\li_n A_n\nnnoojd{C_n}^2A_n^*\Astareta^\mj}^\uujd\negd \nnoduu{X\!\zzmdct{\Dstarti^\emj\negc\sum_\njN\!w_n D_n\nnnoojd{B_n}^2D_n^*\Dstarti^\emj}^\upjd}_2
\lesut \max_{1\lesu n\lesoj\NN}\negc\nnnoou{B_n}\max_{1\lesu n\lesoj\NN}\negc \nnnoou{C_n}\nnnou{X}_2,\label{Cetvorka2}
\end{align}\vsplinc
where  \ineqreff{Cetvorka1} is due to the Cauchy-Schwarz norm inequality (9) in \cite[Lemma~2.1(d1)]{JKL}, applied to the families $\zzvju{\li_n^\ujd\Astareta^\mj A_nC_n}_\njN$ and $\zzvju{w_n^\ejd\Dstarti^\emj D_nB_n}_\njN$ instead of $\AAnN$ and $\BBnN$ therein. The first inequality in \eqref{Cetvorka2} is based on inequalities
$\BnBnz\lesu\nnnoojd{B_n^*}^2I=\nnnoojd{B_n}^2I$ and $\CnCnz\lesu\nnnoojd{C_n^*}^2I=\nnnoojd{C_n}^2I,$
while the second inequality in \eqref{Cetvorka2} is based on the estimates
$\Astareta^\mj\!\sum_\njN\!\li_n\AnAnzu\Astareta^\mj\lesu I$ and $\Dstarti^\emj\sum_\njN\!w_n D_nD_n^*\Dstarti^\emj\lesu I,$ all combined again with the double monotonicity property \eqref{duplamonotonost}.

To prove \ineqreff{jenABCD} we rely on the estimates
\vspbbald\notag
&\!\!\!\!\nnodu{\sum_\njN\!w_n^\ejd A_nC_n\Cdireta^\emj X\Dstarti^\emj D_nB_n}_1
\lesup\sum_\njN\!w_n^\ejd\nnnoju{A_nC_n\Cdireta^\emj X\Dstarti^\emj D_nB_n}_1\\
&\!\!\!\!\lesot\max_{1\lesu n\lesoj\NN}\negc\nnnoou{A_n}\nnnoou{B_n}
\sum_\njN\nnnoju{C_n\Cdireta^\emj X\Dstarti^\emj w_n^\ejd D_n}_1\label{nulaABCDdok}\\
&\!\!\!\!\lesot\max_{1\lesu n\lesoj\NN}\negc\nnnoou{A_n}\nnnoou{B_n}
\nnodtj{\zzmdct{\Cdireta^\emj\!\sum_\njN\!\CnzCn\Cdireta^\emj}^\upjd\negp X\!
\zzmdcd{\Dstarti^\emj\!\sum_\njN\!w_n\DnDnz\Dstarti^\emj}^\upjd}_1
\lespt\max_{1\lesu n\lesoj\NN}\negc\nnnoou{A_n}\nnnoou{B_n}\nnnoou{X}_\ISj\negd,\label{jenABCDdok}
\end{align}\vsplinc
where both inequalities in \eqref{nulaABCDdok} are obvious, while
the first inequality in  \eqref{jenABCDdok} is based on (11) in \cite[Lemma~2.1(e)]{JKL}, applied to the families $\zzvj{C_n\Cdireta^\emj}_{n=1}^\NN$ and $\zzvj{\Dstarti^\emj\pozd w_n^\ejd D_n}_{n=1}^\NN$ instead of $\AAnN$ and $\BBnN.$ The last inequality in \eqref{jenABCDdok} is based on inequalities
$\Cdireta^\emj\!\sum_\njN\!\CnzCn\Cdireta^\emj\lesu\Cdireta^\emj\Cdireta^2\Cdireta^\emj\equ I$ and (similarly derived) $\Dstarti^\emj\!\sum_\njN\!w_n\DnDnz\Dstarti^\emj\lesu I,$ both additionally combined with the double monotonicity property \eqref{duplamonotonost}.

To prove that \ineqreff{nulaABCD} implies \ineqreff{nulaABCDalt}, we first realize that for any $\eta,\zeta>0$
\vspbbald\notag
\nnodu{\sum_\njN\!\li_n^\ujd\AnCnXDnBn}
&=\nnodu{\Astareta\Astareta^\mj\negc\sum_\njN\!\li_n^\uujd\AnCnXDnBn\Bdirzi^\emj\Bdirzi}\\
&\lespj\nnnoou{\Astareta}\nnodu{\Astareta^\mj\negc\sum_\njN\!\li_n^\uujd\AnCnXDnBn\Bdirzi^\emj}\!\nnnoou{\Bdirzi}\notag\\
&\lespj\max_{1\lesu n\lesoj\NN}\negc\nnnoou{C_n}\nnnoou{D_n}
\nnodu{\eta I\pludp\sum_\njN\!\li_n\AnAnz}^\uujd\negd\nnodu{\zeta I\pludp\sum_\njN\!\BnzBn}^\uujd\negd\nnnoou{X},
\label{od2do6}
\end{align}\vsplinc
where the (last) inequality in \eqref{od2do6} is due to \ineqreff{nulaABCD}.
As $\inf_{\eta>0}\nnnojod{\eta I\pludc\sum_\njN\!\li_n\AnAnz}=\nnnojod{\sum_\njN\!\li_n\AnAnz}$ and
$\inf_{\zeta>0}\nnnojod{\zeta I\pludc\sum_\njN\!\BnzBn}=\nnnojod{\sum_\njN\!\BnzBn},$ then by taking $\inf_{\eta>0}$
and $\inf_{\zeta>0}$ on the righthand side of the inequality in \ineqreff{od2do6} we confirm the validity of \ineqreff{nulaABCDalt}.
To prove the opposite implication, an application of  \ineqreff{nulaABCDalt} to $\odjuu{\Astareta^\mj A_n, B_n\Bdirzi^\emj}$ instead of $\oduu{A_n,B_n},$ combined with inequalities $\Astareta^\mj\negc\sum_\njN\!\li_n\AnAnz\Astareta^\mj\les I$ and $\Bdirzi^\emj\negc\sum_\njN\!\BnzBn\Bdirzi^\emj\les I,$ gives
\vspbbald\notag
&\nnodu{\Astareta^\mj\negc\sum_\njN\!\li_n^\uujd\AnCnXDnBn\Bdirzi^\emj}\\
&\lespj\max_{1\lesu n\lesoj\NN}\negc\nnnoou{C_n}\nnnoou{D_n}
\nnodu{\Astareta^\mj\negc\sum_\njN\!\li_n\AnAnz\Astareta^\mj}^\uujd\negd
\nnodu{\Bdirzi^\emj\negc\sum_\njN\!\BnzBn \Bdirzi^\emj}\nnnoou{X}
\lestj \max_{1\lesu n\lesoj\NN}\negc\nnnoou{C_n}\nnnoou{D_n}\nnnoou{X},\notag
\end{align}\vsplintp
validating \ineqreff{nulaABCDalt} and confirming the equivalence of \ineqreff{nulaABCD} and \ineqreff{nulaABCDalt}.

Similarly, starting from \ineqreff{dvaABCD2}, for every $\eta,\zeta>0$ we get
\vspbbald\notag
\nnodu{\sum_\njN\!\AnCnXDnBn}_2
&=\nnodu{\sum_\njN\!A_nC_n\Cdireta^{\pozd\mj}\Cdireta XD_nB_n\Bdirzi^\emj\Bdirzi}_2\\
&\lesoj \nnodu{\sum_\njN\!A_nC_n\Cdireta^{\pozd\mj}\Cdireta XD_nB_n\Bdirzi^\emj}_2\nnnoou{\Bdirzi}
\lesoj \max_{1\lesu n\lesoj\NN}\negc\nnnoou{A_n}\max_{1\lesu n\lesoj\NN}\negc\nnnoou{D_n}\nnnoou{\Cdireta X}_2\nnnoou{\Bdirzi}\notag\\
&=\max_{1\lesu n\lesoj\NN}\negc\nnnoou{A_n}\max_{1\lesu n\lesoj\NN}\negc\nnnoou{D_n}
\nnodjd{\zeta I\plusdp\sum_\njN\!\BnzBn}^\uujd\nnodtj{\zzmduu{\eta I\plusdp\sum_\njN\!\CnzCn}^\upjd\negp X}_2,\notag
\end{align}\vsplint
so by realizing that
$\inf_{\zeta>0}\nnnojod{\zeta I\plusdc\sum_\njN\!\BnzBn}=\nnnojod{\sum_\njN\!\BnzBn}$ and
$\inf_{\eta>0}\nnnojtd{\zzmju{\eta I\plusdp\sum_\njN\!\CnzCn}^\uujd\negd X}_2
=\nnnojtj{\zzmju{\sum_\njN\!\CnzCn}^\uujd\negd  X}_2,$
 this confirms that  \ineqreff{dvaABCD2} implies \ineqreff{dvaABCD2alt}. Indeed, the last equality follows from   $\nnnojtd{\zzmj{\eta_1 I+\sum_\njN\!\CnzCn}^\uujd X}_2\les\nnnojtd{\zzmj{\eta_2 I\plusdc\sum_\njN\!\CnzCn}^\uujd\!X}_2$ for $0<\eta_1<\eta_2,$ which follows from the double monotonicity property \eqref{duplamonotonost} in the Hilbert-Schmidt ideal $\ccd,$ applied on $A\dodu
\zzmju{\eta_1 I+\sum_\njN\!\CnzCn}^\uujd,$ $C\dodu\zzmjuu{\eta_2 I\pludc\sum_\njN\!\CnzCn}^\uujd$ and $B\dodu D\dod I.$ Also,
\vspbbald\notag
&\nnodtj{\zzmddt{\eta I\plusdp\sum_\njN\!\CnzCn}^\upjd\negp X}_2 -\nnodtj{\zzmddt{\sum_\njN\!\CnzCn}^\upjd\negp X}_2
\lesu\nnodtj{\zzmddt{\eta I\plusdp\sum_\njN\!\CnzCn}^\upjd\negp X\mintt\zzmddt{\sum_\njN\!\CnzCn}^\upjd\negp X}_2\\
&\lesu\nnodtj{\zzmddt{\eta I\plusdp\sum_\njN\!\CnzCn}^\upjd \minpd\zzmddt{\sum_\njN\!\CnzCn}^\upjd}\nnnoo{X}_2
\lesu\nnodjd{\eta I\plusdp\sum_\njN\!\CnzCn\mindp\sum_\njN\!\CnzCn}^\uujd \negd\nnnoo{X}_2
=\kor{\eta}\nnnoo{X}_2,\label{5do12}
\end{align}\vsplintp
where the last inequality in \eqref{5do12} is due to \cite[Cor.~2(i)]{A}, implying that
the function $f\colon\Rp\to\Rp\colon\eta\mapsto\allowbreak\nnnojtd{\zzmj{\eta I\pludc\sum_\njN\!\CnzCn}^\uujd\!X}_2$ is continuous and increasing on $\Rp,$ implying $\inf_{\eta>0}f(\eta)=\lim_{\eta\ssearrow0}f(\eta)=f(0)$ and
justifying \eqref{5do12}.

For the opposite implication, by using \ineqreff{dvaABCD2alt} for $\Cn\Cdireta^{\pozd\mj}$ and $B_n\Bdirzi^\emj$
instead of $C_n$ and $B_n,$ combined with the inequalities $\Bdirzi^{\pozd\mj}\negc\sum_\njN\!\BnzBn\Bdirzi^\emj\lesu I$ and $\Cdireta^{\pozd\mj}\negc\sum_\njN\!\CnzCn\Cdireta^\emj\lesu I,$ we derive
\vspbbald\notag
&\nnodu{\sum_\njN\!A_n^\fff\Cn\Cdireta^\emj X\Dn\Bn\Bdirzi^\mj}_2\notag\\
&\lespj\max_{1\lesu n\lesoj\NN}\negc\nnnoou{A_n}\negd\max_{1\lesu n\lesoj\NN}\negc \nnnoou{D_n}\nnodu{\Bdirzi^{\pozd\mj}\negc\sum_\njN\!\BnzBn\Bdirzi^{\pozd\mj}}^\uujd
\nnodcd{\zzmdtd{\Cdireta^\emj\negc\sum_\njN\!\CnzCn\Cdireta^\emj}^\upjd\negc X}_2
\lespj\max_{1\lesu n\lesoj\NN}\negc\nnnoou{A_n}\negd\max_{1\lesu n\lesoj\NN}\negc \nnnoou{D_n}\nnnoou{X}_\ISd.\notag
\end{align}\vsplinc
So \ineqreff{dvaABCD2alt} implies \ineqreff{dvaABCD2}.

The equivalence of \ineqreff{dvaABCD} and \ineqreff{dvaABCDalt} proves similarly, with the same
arguments used in the proof of the equivalence of \ineqreff{dvaABCD2} and \ineqreff{dvaABCD2alt}.

To prove that \ineqreff{jenABCD} implies \ineqreff{jenABCDalt}, note that for any $\eta,\ti >0$
\vspbbald\notag
&\nnodjd{\sum_\njN\!w_n^\ejd A_n^\fff\Cn X\Dn\Bn}_1
=\nnodjd{\sum_\njN\!w_n^\ejd A_n^\fff\Cn \Cdireta^\emj\Cdireta X\Dstarti\Dstarti^\emj\Dn\Bn}_1\\
&\lespj\max_{1\lesu n\lesoj\NN}\negc\nnnoojd{A_n^\fff}\nnnoojd{\Bn}\nnnoou{\Cdireta X\Dstarti}_\ISj
=\max_{1\lesu n\lesoj\NN}\negc\nnnoojd{\Anff}\nnnoojd{\Bn}\nnodtd{\zzmdjc{\eta I\pludc\sum_\njN\!\CnzCn}^\upjd\negp X\! \zzmdjc{\ti I\pludc\sum_\njN\!w_n\DnDnz}^\upjd}_1,\label{jenABCDaltDOK}
\end{align}\vsplinc
where \ineqreff{jenABCDaltDOK} follows from \ineqreff{jenABCD} applied to $\Cdireta X\Dstarti$ instead of $X.$

It was shown in \eqref{5do12} that
$\nnnojtj{\zzmjjj{\eta I\pludc\sum_\njN\!\CnzCn}^\utjd\minpj\zzmjjj{\sum_\njN\!\CnzCn}^\upjd}
\lesu\kor{\eta}\to 0$ as $\eta\ssearrow 0,$ and similarly
$\nnnojtj{\zzmjjj{\ti I\pludc\sum_\njN\!w_n\DnDnz}^\utjd\negp-\zzmjjj{\sum_\njN\!w_n\DnDnz}^\utjd}
\lesu\kor{\ti}\to 0$ as $\ti\ssearrow0,$ which therefore easily implies that
$\lim\limits_{\eta,\ti\ssearrow0}\negd\nnnojtj{\zzmju{\eta I\pludc\sum_\njN\!\CnzCn}^\utjd\negd X\negd
\zzmju{\ti I\pludc\sum_\njN\!w_n\DnDnz}^\utjd}_\ISj\negd
=\nnnojtj{\zzmju{\sum_\njN\!\CnzCn}^\utjd X \zzmju{\sum_\njN\!w_n\DnDnz}^\utjd}_\ISj.$
So by taking the limit as $\eta,\ti\ssearrow0$ on the righthand side of \eqref{jenABCDaltDOK},
\ineqreff{jenABCDalt} follows.

To prove that $\eqref{jenABCDalt}$ implies $\eqref{jenABCD},$ by applying \eqref{jenABCDalt} to $\odju{C_n\Cdireta^\emj,\Dstarti^\emj D_n}$ instead of $\odu{C_n,D_n}$ we get
\vspbbald\notag
&\nnodu{\sum_\njN\!w_n^\ejd A_n^\fff\Cn\Cdireta^\emj X\Dstarti^\mj\Dn\Bn}_1\\
&\lespj\max_{1\lesu n\lesoj\NN}\negc\nnnoou{A_n^\fff}\nnnoou{\Bn} \nnodu{\zzmju{\Cdireta^\emj\sum_\njN\!\CnzCn\Cdireta^\emj}^\uujd X\zzmju{\Dstarti^\emj\sum_\njN\DnDnz\Dstarti^\emj}^\uujd}_1
\lespj\max_{1\lesu n\lesoj\NN}\negc\nnnoou{A_n^\fff}\nnnoou{\Bn}\nnnoou{X}_1\notag,
\end{align}\vsplinc
since $\Cdireta^\emj\sum_\njN\!\CnzCn\Cdireta^\emj$ and $\Dstarti^\emj\sum_\njN w_n\DnDnz\Dstarti^\emj \lesu I,$ again combined with the property \eqref{duplamonotonost}.

To prove that \eqref{nulaABCD} implies \eqref{nula}, we just take $C_n\dodu D_n\dodu I$ in \eqref{nulaABCD} for all $\eN\ni n\lesu\N,$ while by taking $A_n\dodu D_n\dodu I,$  $C_n\dodu A_n$ and $\eta\dodu\ei$ in \eqref{dvaABCD2} for any $\eN\ni n\lesu\N$ we obtain  \ineqreff{dva2}.
Similarly, by taking $B_n\dodu C_n\dodu I$ and $D_n\dodu B_n$ (resp.\ $A_n\dodu B_n\dodu I, C_n\dodu A_n, D_n\dodu B_n$ and $\eta\dod\ei$) for any $\eN\ni n\lesu\N$ in \eqref{dvaABCD} (resp. \eqref{jenABCD})
we get  \ineqreff{dva} (resp.\ \ineqreff{jen}).
\eeproo\medskip

\section{Main results}

\bbt\label{CSinCp}
Let $q,r,s\gesuu1\!$ satisfy $\!\frac1{2q}\plu\frac1{2r}\equt\frac1s\!$ and let also $\!\nizjd{\lambda},\!\nizjd{w}\!\!$ be sequences in $\!\Rpz.\!\!$ If $\nizjd{A},\odjd{\li_n^\utjd\!\Anuz}_\nji,\nizjd{B},\odjd{w_n^\uujd\!\Bnuz}_\nji\inu\ldsBH$ and $X\inu\ccs,$ then
$\TSCssummuc\li_n^{\ucfrac12\mindj\frac1{2q}}\!w_n^\frac1{2r}\!\AnXBn$ exists and
\vspbbald\label{weghtedCSforCp}\!
&\sup_{\eta,\zi\negd\IIS>\negj0}\negj
\nnodct{\zzmdct{\eta I\plutp\jsummup\li_n\AnAnz}^{\uvifrac1{2q}\minu\frac12}\negp
\Cssumm\!\!\li_n^{\ucfrac12\mindj\frac1{2q}}\!w_n^\uufrac1{2r}\negp\AnXBn\!
\zzmduc{\zi I\plutp\jsummup\BnzBn}^{\uvifrac1{2r}\minu\frac12}}_{\negd s}
\lesst\nnodct{\zzmddc{\jsummup\AnzAn}^\uvifrac1{2q}\negp X\!\zzmddc{\jsummup w_n\BnBnz}^\uvifrac1{2r}}_{\negd s}\negc,\\
&\nnodcd{\Cssumm\!\!\li_n^{\ucfrac12\mindj\frac1{2q}}\!w_n^\frac1{2r}\!\AnXBn}_{\negd s}
\lessd\nnodod{\jsummuc\li_n\AnAnz}^{\uufrac12\minu\frac1{2q}}\!
\nnodjd{\jsummup\BnzBn}^{\uufrac12\minu\frac1{2r}}\!
\nnodct{\zzmddt{\jsummup\AnzAn}^\uvifrac1{2q}\negp X\!
\zzmddc{\jsummup w_n\BnBnz}^\uvifrac1{2r}}_{\negd s}\!.\label{weghtedCSforCpAlt}
\end{align}\vsplinc
\eet
\bbproo
Based on the notation from Definition \ref{defInvertibilni}, we consider the function
\vspbbaldp\label{biStripFtion}
\vfi\colon\bbS\OJ\xuu\bbS\OJ\to\BH\colon
\oduu{z,v}\mto\sum_\njN\!\li_n^\frac{1\mindj z}2\! w_n^\frac{v}2 \Astareta^{z\minjd1}\!A_n\Adirei^{-z} X\Bstarti^{-v}\Bn\Bdirzi^{v\minu1}\!,
\end{align}\vsplincp
which is continuous on its domain and holomorphic in each variables on its interior $\bbS\OJo\xuu\bbS\OJo,$ 
and which for all $s,t\in\eR$ satisfies the following boundary conditions:
\vspbbald\label{nulaABocena}
\nnnoojd{\vfi\odu{is,it}}
&\eqdj\nnodjd{\sum_\njN\!\li_n^\uvifrac{1\mindj is}2
\!w_n^\utfrac{it}2\! \Astareta^{is\minjd1}\!A_n\Adirei^{\minu is}X
\Bstarti^{\minu it}\!B_n\Bdirzi^{it\minu1}}\lesdj\nnnoojd{X},\\
\nnnoou{\vfi\odu{1\pludj is,it}}_\ISd&\eqcj\nnodjd{\sum_\njN\!\li_n^\uvifrac{\minu is}2\negc w_n^\utfrac{it}2\! \Astareta^{is} A_n\Adirei^{\minuu1\minjd is}\!X
\Bstarti^{\minu it}\!B_n\Bdirzi^{it\minjd1}}_\ISd
\lespj\nnnoou{X}_\ISd,\label{dvaABocena}\\
\nnnoou{\vfi\odu{is,1\pludj it}}_\ISd&\eqcj\nnodjd{\sum_\njN\!\li_n^\uvifrac{1\mindj is}2\!
w_n^\utfrac{1\pludj it}2\!\! \Astareta^{is\minjd1}\!A_n\Adirei^{\minu is}\!X
\Bstarti^{\minuu1\minu it}\negp\Bn\Bdirzi^{it}}_\ISd
\lespj\nnnoou{X}_\ISd,\label{altdvaABocena}\\
\nnnoou{\vfi\odu{1\pludj is,1\plusdj it}}_\ISj&\eqcj\nnodjd{\sum_\njN\!\li_n^\uvifrac{\minu is}2\negc
w_n^\utfrac{1\pludj it}2\! \Astareta^{is} A_n\Adirei^{\minuu1\mindj is}\!X
\Bstarti^{\minuu1\minjd it}\!B_n\Bdirzi^{it}}_\ISj
\lespj\nnnoou{X}_\ISj.\label{jenABocena}
\end{align}\vsplind
\ineqreff{nulaABocena} (resp.\ \ineqreff{dvaABocena}) follows from  \ineqreff{nula}
(resp.\ \ineqreff{dva2}) applied to $(\li_n^\uvifrac{\minu is}2\negp\Astareta^{is}\negj A_n,
w_n^\utfrac{it}2\! B_n \Bdirzi^{it},\allowbreak\Adirei^{\minu is}X
\Bstarti^{\minu it})$
(resp.\ $(\li_n^\uvifrac{\minu is}2\!\Astareta^{is} A_n,w_n^\utfrac{it}2\!B_n\Bdirzi^{it},
\Adirei^{\minu is}\!X\Bstarti^{\minu it})$) instead of $(A_n,B_n,X),$ supported by the fact that
$\odjjc{\eta I\pludc\sum_\njN\negc\li_n\li_n^\uvifrac{\minu is}2\!\!\Astareta^{is}\negj A_n\oddt{\li_n^\uvifrac{\minu is}2\!\!\Astareta^{is}\negd A_n}^{\!*}}^\ucfjd
\equc\odju{\eta I\plutc\sum_\njN\negc\li_n\Astareta^{is}\negj\AnAnz\Astareta^{\minjd is}}^\ucfjd
\equc\odjuu{\Astareta^{is}\negd\zzmjjj{\eta I\plutc\sum_\njN\negc\li_n\AnAnz}\Astareta^{-is}}^\ucfjd \negc\allowbreak
=\negc\zzmj{\Astareta^{is}\Astareta^2\Astareta^{-is}}^\ucfjd\!=\Astareta,$
$\zzmj{\zi I\pludc\sum_\njN\! \zzmjdj{w_n^\utfrac{it}2\!
B_n \Bdirzi^{it}}^{\!*}\!w_n^\utfrac{it}2\! B_n \Bdirzi^{it}}^\ucfjd\!=\negd\Bdirzi$
and also $\nnnoou{\Adirei^{\minu is}\negj X\Bstarti^{\minu it}}\negd =\negd\nnnoou{X}$ (resp.\
$\zzmj{\ei I\plujj\sum_\njN\!\zzmotd{\li_n^\uvifrac{\minu is}2\!\Astareta^{is}\negj A_n}^{\!*} \negd\li_n^\uvifrac{\minu is}2\!\Astareta^{is}\negj A_n}^\ucfjd=\Adirei,$
$\zzmj{\zi I\plujj\sum_\njN\!
\zzmjdj{w_n^\utfrac{it}2\!B_n\Bdirzi^{it}}^{\!*}\!w_n^\utfrac{it}2\!B_n\Bdirzi^{it}}^{\ucfjd}=\Bdirzi^{}$ and also $\nnnoou{\Adirei^{\minu is}\!X\Bstarti^{\minu it}}_\ISd=\nnnoou{X}_\ISd$).
Similarly, by applying \ineqreff{dva} (resp.\ \ineqreff{jen}) on
$(\li_n^\uvifrac{\mindj is}2\negp\Astareta^{is}\negj A_n,
w_n^\utfrac{it}2\Bstarti^{\minu it}\Bn\Bdirzi^{it},\allowbreak
\Adirei^{\minu is}\!X)$ (resp.\ on
$(\negj\li_n^\uvifrac{\minu is}2\negc\Astareta^{is} A_n,
w_n^\utfrac{it}2\!B_n\Bdirzi^{it},\Adirei^{\mindj is}\!X
\Bstarti^{\minjd it})$), instead on $(A_n,B_n,X),$ we obtain \ineqreff{altdvaABocena} (resp.\ \ineqreff{jenABocena}).

For any $\ai,\bi\inu\OJ$ by $[\ai]$ interpolating \ineqreff{nulaABocena} and \ineqreff{dvaABocena} it follows
\vspbbald\notag
&\negxiv\nnnoojd{\vfi\odu{\ai\plu is,it}}_{_{{\IS p}_\ai}}\eqpd
\nnnoojd{\vfi\odu{is\odu{1\mindj\ai}\plu\odu{1\plu is}\ai,it}}_{_{{\IS p}_\ai}}\\
&\negxiv\eqpd\nnodjt{\sum_\njN\!\li_n^\uvifrac{\odu{1\mindj is} \odu{1\mindj\ai}-is\ai}2\!w_n^\utfrac{it}2\!\Astareta^{\odu{is\minjd1}\odu{1\mindj\ai}\plu is\ai}\negp \Anff\Adirei^{\minu is\odu{1\mindj\ai}\minu\odu{1\pludj is}\ai}\!X
\Bstarti^{\minu it}\negd\Bn\Bdirzi^{it\minu1}}_{p_\ai}\notag\\
&\negxiv\eqpd\nnodjt{\sum_\njN\!\li_n^\uvifrac{1\mindj\ai-is}2\! w_n^\utfrac{it}2\!
\Astareta^{\ai\plu is\minjd1}\negp\Anff\Adirei^{\minu is\minjd1}\!X
\Bstarti^{\minu it}\negd\Bn\Bdirzi^{it\minu1}}_{p_\ai}
\lespd\nnnoojd{X}_{_{{\IS p}_\ai}}\!,
\label{interpolocena1}
\end{align}\vsplinc
while by $[\ai]$ interpolating \ineqreff{altdvaABocena} and \ineqreff{jenABocena} we get
\vspbbald\notag
\negxiv\nnnoojd{\vfi\odu{\ai\plu is,1\plu it}}_{_{{\IS q}_\ai}}
&\eqsd\nnodjd{\sum_\njN\!\li_n^\uvifrac{\odu{1\mindj is}\odu{1\mindj\ai}-is\ai}2\!
w_n^\frac{1\pludj it}2\!\!\Astareta^{\odu{is\minjd1}\odu{1\mindj\ai}\plu is\ai}\negp\Anff
\Adirei^{\minu is\odu{1\mindj\ai}\minu\odu{1\pludj is}\ai}\!X
\Bstarti^{\minjd1\minu it}\!\Bn\Bdirzi^{it}}_{q_\ai}\notag\\
&\eqsd\nnodjd{\sum_\njN\!\li_n^\uvifrac{1\mindj\ai-is}2\!
w_n^\frac{1\pludj it}2\!\!\Astareta^{\ai\plu is\minjd1}\negp\Anff\Adirei^{\minu\ai\minu is}\!X
\Bstarti^{\minjd1\minu it}\!\Bn\Bdirzi^{it}}_{q_\ai}
\lessd\nnnoou{X}_{_{{\IS q}_\ai}}\!,\!\label{interpolocena2}
\end{align}\vsplinc
where $\frac1{p_\ai}\dodu\odu{1\mindj\ai}\cdot0+\ai\cdot\frac12=\frac{\ai}2$ and
$\frac1{q_\ai}\dodu\odu{1\mindj\ai}\cdot\frac12+\ai\cdot1=\frac{1\pludj\ai}2.$
Finally, by $[\bi]$ interpolating \ineqreff{interpolocena1} and \ineqreff{interpolocena2} we obtain
\vspbbald\notag
\!\!\!\!\nnnoojd{\vfi\odu{\ai\plu is,\bi\plu it}}_{_{{\IS r}_{\negp_\bi}}}
\!\!&\eqsd\nnodjd{\sum_\njN\!\li_n^\uvifrac{1\mindj\ai\minu is}2\!
w_n^{\odu{1\mindj\bi}\frac{it}2\plu\bi\frac{1\pludj it}2}\!\!\Astareta^{\ai\minjd1\plu is}\negc\Anff
\Adirei^{\minu\ai\minu is}\!X
\Bstarti^{\minu\odu{1\mindj\bi}{it}\minu\bi\odu{1\pludj it}}\!\Bn
\Bdirzi^{\odu{1\mindj\bi}\odu{it\minjd1}\plu\bi{it}}}_{r_{\negp_\bi}}\\
&\eqsd\!\nnodjd{\sum_\njN\!\li_n^\uvifrac{1\mindj\ai\minu is}2\!
w_n^\frac{\bi\plu it}2\!\!\Astareta^{\ai\minjd1\plu is}\negc\Anff\Adirei^{\minu\ai\minu is}\!X
\Bstarti^{\minu\bi\minu it}\!\Bn\Bdirzi^{\bi\plu it\minjd1}}_{r_{\negp_\bi}}
\lessd\nnnoou{X}_{_{{\IS r}_{\negp_\bi}}}\!,\!\label{interpolocena3}
\end{align}\vsplint
where $\frac1{r_{\negp_\bi}}\dodu\odu{1\mindj\bi}\cdot\frac1{p_\ai}\plu\bi\cdot\frac1{q_\ai}
=\frac{\odu{1\mindj\bi}\ai}2\plu\frac{\bi\odu{1\plu\ai}}2=\frac{\ai\plu\bi}2.$
By taking the special case $s\dodu t\dodu0$ and using the change of variable
$X\mto\Adirei^{\ai\plu is}\!X\Bstarti^{\bi\plu it}$ 
in the inequality in \eqref{interpolocena3} it implies
\vspbbald
\nnodjd{\sum_\njN\!\li_n^\uvifrac{1\mindj\ai}2\!
w_n^\frac{\bi}2\!\Astareta^{\ai\minjd1}\!A_nXB_n\Bdirzi^{\bi\minjd1}}_{r_{\negp_\bi}}
\lessd\nnnoju{\Adirei^{\ai\plu is}\!X\Bstarti^{\bi\plu it}}_{r_{\negp_\bi}}
\eqsd\nnnoju{\Adirei^\ai X\Bstarti^\bi}_{r_{\negp_\bi}}\!,\!\label{interpolocena4}
\end{align}\vsplinc
so by letting $\ei\sssearrow0,\ti\sssearrow0$ and by using the additional change of variables
$\odu{q,r,s}\dodu\odju{1\kroz\ai,1\kroz\bi,r_{\negp_\bi}},$ then \ineqreff{interpolocena4} can be rewritten as
\vspbbald
\!\!\!\!\nnodut{\zzmdct{\eta I\pludp\sum_\njN\!\li_n\AnAnz}^{\negp\frac1{2q}\minu\frac12}\negp
\sum_\njN\!\li_n^{\frac12\mindj\frac1{2q}}\!w_n^\frac1{2r}\!A_nXB_n\!
\zzmddc{\zi I\pludp\sum_\njN\!\BnzBn}^{\negp\frac1{2r}\minu\frac12}}_s
\lessd\nnodut{\zzmddt{\sum_\njN\!\AnzAn}^{\negp\frac1{2q}}\negp X\!\zzmddc{\sum_\njN\!w_n\BnBnz}^{\negp\frac1{2r}}}_s\!,\!\!\label{ocena5}
\end{align}\vsplinc
which implies
\vspbbald
\!\!\!\!\nnoduu{
\sum_\njN\!\li_n^{\frac12\mindj\frac1{2q}}\!w_n^\frac1{2r}\!A_nXB_n}_s
\lessd\nnoduu{\sum_\njN\!\li_n\AnAnz}^{\!\frac12\minu\frac1{2q}}\!
\nnoduu{\sum_\njN\!\BnzBn}^{\!\frac12\minu\frac1{2r}}\!
\nnodut{\zzmddt{\sum_\njN\!\AnzAn}^{\negp\frac1{2q}}\negp X\!\zzmddc{\sum_\njN\!w_n\BnBnz}^{\negp\frac1{2r}}}_s\!.\!\!\label{ocena6}
\end{align}\vsplinc
Indeed, we have
\vspbbald
&\nnoduu{\sum_\njN\!\li_n^{\utfrac12\mindj\frac1{2q}}\!w_n^\frac1{2r}\!A_nXB_n}_s
\lespd\nnodcd{\zzmdct{\eta I\pludp\sum_\njN\!\li_n\AnAnz}^{\uvifrac12\minu\frac1{2q}}}\negj
\nnodcd{\zzmduc{\zi I\pludp\sum_\njN\!\BnzBn}^{\uvifrac12\minu\frac1{2r}}}\notag\\
&\,\,\,\x\nnodut{\zzmdct{\eta I\pludp\sum_\njN\!\li_n\AnAnz}^{\uvifrac1{2q}\minu\frac12}\negp
\sum_\njN\!\li_n^{\utfrac12\mindj\frac1{2q}}\!w_n^\frac1{2r}\!A_nXB_n
\zzmduc{\zi I\pludp\sum_\njN\!\BnzBn}^{\uvifrac1{2r}\minu\frac12}}_s\label{njn1}\\
&\eqod\nnodjd{\eta I\pludp\sum_\njN\!\li_n\AnAnz}^{\uufrac12\minu\frac1{2q}}\!
\nnodjd{\zi I\pludp\sum_\njN\!\BnzBn}^{\uufrac12\minu\frac1{2r}}\notag\\&
\,\,\,\x\nnodut{\zzmdct{\eta I\pludp\sum_\njN\!\li_n\AnAnz}^{\uvifrac1{2q}\minu\frac12}\!\! \sum_\njN\!\li_n^{\utfrac12\mindj\frac1{2q}}\!w_n^\frac1{2r}\!A_nXB_n
\zzmduc{\zi I\pludp\sum_\njN\!\BnzBn}^{\negp\frac1{2r}\minu\frac12}}_s\label{njn2}\\
&\lesod\nnodjd{\eta I\pludp\sum_\njN\!\li_n\AnAnz}^{\uufrac12\minu\frac1{2q}}\!
\nnodjd{\zi I\pludp\sum_\njN\!\BnzBn}^{\uufrac12\minu\frac1{2r}} \nnodut{\zzmddt{\sum_\njN\!\AnzAn}^{\uvifrac1{2q}}\negp X\!\zzmddc{\sum_\njN\!w_n\BnBnz}^{\uvifrac1{2r}}}_s.\label{njn3}
\end{align}\vsplinc
Here \ineqreff{njn1} follows from property $\nnnoo{CY\negd D}_s\lesu\nnnoo{C}\nnnoo{Y}_s\nnnoo{D}$ for $C,D\in\BH$ and $Y\inu\ccs,$ \eqreff{njn2} is based on property $\nnnoou{A^p}=\nnnoo{A}^p$ for positive operator $A$ and $p>0,$ which directly follows from the spectral mapping theorem $f(\sigma(A))=\sigma(f(A))$ for continuous function $f$ combined with the $\|A\|=\max_{\lambda\inu\sigma(A)}\lambda$, while \ineqreff{njn3} follows from \ineqreff{ocena5}.
So by letting $\eta\sssearrow0,\zi\sssearrow0$ in \eqref{njn3} we obtain \ineqreff{ocena6}.

Based on the previous inequalities \eqref{njn1}-\eqref{njn3}, for any $\en\ni\M>\N$
we get the first inequality in
\vspbbald \notag
&\nnoduu{\sum_{k=\NN\pluu1}^\MM\!\!\li_k^{\frac12\mindj\frac1{2q}}\!w_k^\frac1{2r}\!A_kXB_k}_s
\lespd\nnodjd{\eta I\pludp\negx\sum_{k=\NN\pluu1}^\MM\!\!\!\li_k\AkAkz}^{\utfrac12\minu\frac1{2q}}\!
\nnodjt{\zi\negd I\pludp\negx\sum_{k=\NN\pluu1}^\MM\!\!\BkzBk}^{\utfrac12\minu\frac1{2r}}\!
\nnodut{\zzmddt{\sum_{k=\NN\pluu1}^\MM\!\!\!\AkzAk}^{\negp\frac1{2q}}\negp X\!\zzmddc{\sum_{k=\NN\pluu1}^\MM\!\!\!w_k\BkBkz}^{\negp\frac1{2r}}}_s\\
&\lesod\nnodjd{\eta I\pludp\jsummup\li_n\AnAnz}^{\frac12\minu\frac1{2q}}\!
\nnodjt{\zi I\pludp\jsummup\BnzBn}^{\uufrac12\minu\frac1{2r}}\!
\nnodut{\zzmddt{\jsummNup\AkzAk}^\uvifrac1{2q}\negp X\!
\zzmddc{\jsummNup w_k\BkBkz}^\uvifrac1{2r}}_s\to 0\label{njn4}
\end{align}\vsplinc
as $\N\to\iip,$ since $\jlim{\!\NN\to\iip}\!\!\TSjsummNuc\AkzAk=\jlim{\!\NN\to\iip}\!\!\TSjsummNuc w_k\BkBkz\equu0$ and $X\inuu\ccs.$
The second inequality in \eqref{njn4} is based on the monotonicity of operator norm $\nnnou{\cdot}$
for positive operators, as well as on the double monotonicity property \eqref{duplamonotonost} for $\nnnou{\cdot}_s$ norm, combined with the operator monotonicity of the functions $t\mto t^\frac1q$ and $t\mto t^\frac1r$ on $\Rp;$ see \cite[Ex.~2.5.9(3)]{Hi}. As $\ccs$ is a Banach space, this implies that
there exists $\TSCssumm\negc\li_n^{\ucfrac12\mindj\frac1{2q}}\!w_n^\frac1{2r}\!\AnXBn.$

The proof of \eqref{weghtedCSforCpAlt} is based on following chain of inequalities:
\vspbbaldp
&\nnodcd{\Cssumm\negc\li_n^{\utfrac12\mindj\frac1{2q}}\!w_n^\frac1{2r}\negc A_nXB_n}_s
\eqpd\nnodsd{\Cslimm\sum_\njN\!\li_n^{\utfrac12\mindj\frac1{2q}}\!w_n^\frac1{2r}\negc A_nXB_n}_s
\equp\lim_{\NN\to\iii}\negd\nnodjd{\sum_\njN\!\li_n^{\utfrac12\mindj\frac1{2q}}\!w_n^\frac1{2r}\!A_nXB_n}_s\notag\\
&\lesoc\lim_{\NN\to\iii}\nnoduu{\sum_\njN\!\li_n\AnAnz}^{\uufrac12\minu\frac1{2q}}\!
\nnoduu{\sum_\njN\!\BnzBn}^{\uufrac12\minu\frac1{2r}}\!
\nnodut{\zzmddt{\sum_\njN\!\AnzAn}^\uvifrac1{2q}\negp X\! \zzmddc{\sum_\njN\!w_n\BnBnz}^\uvifrac1{2r}}_s\notag\\
&\lesoj\nnodjd{\jsummup\!\li_n\AnAnz}^{\uufrac12\minu\frac1{2q}}\!
\nnodjd{\jsummup\!\BnzBn}^{\uufrac12\minu\frac1{2r}}\!
\nnodut{\zzmddt{\jsummup\!\AnzAn}^\uvifrac1{2q}\negp X\!\zzmddc{\jsummup\!w_n\BnBnz}^\uvifrac1{2r}}_s,  \label{Finale}
\end{align}
where the first inequality in \eqref{Finale} is according to \eqref{ocena6}, while the last inequality in \eqref{Finale} is based on the same arguments used in the proof of the last inequality in \eqref{njn4}.

Now, for all $\eta,\zi>0$
\vspbbald\notag
&\nnodct{\zzmdct{\eta I\pludp\jsummup\li_n\AnAnz}^{\uvifrac1{2q}\minu\frac12}\negp
\Cssummuc\li_n^{\ucfrac12\mindj\frac1{2q}}\!w_n^\uufrac1{2r}\!\AnXBn\!
\zzmdtc{\zi I\pludp\jsummup\BnzBn}^{\uvifrac1{2r}\minu\frac12}}_s\\\label{finis1}
&\lesod\nnodct{\zzmdct{\eta I\pludp\jsummup\li_n\AnAnz}^{\uvifrac1{2q}\minu\frac12}}\negd
\nnodct{\zzmdtc{\zi I\pludp\jsummup\BnzBn}^{\utfrac1{2r}\minu\frac12}}\negd \nnodcd{\Cssummuc\li_n^{\ucfrac12\mindj\frac1{2q}}\!w_n^\uufrac1{2r}\!\AnXBn}_s\\\label{finis2}
&\eqod\nnodjt{\eta I\pludp\jsummup\li_n\AnAnz}^{\utfrac1{2q}\minu\frac12}\negc
\nnodjt{\zi I\pludp\jsummup\BnzBn}^{\utfrac1{2r}\minu\frac12}\negc
\nnodcd{\Cssummuc\li_n^{\ucfrac12\mindj\frac1{2q}}\!w_n^\uufrac1{2r}\!\AnXBn}_s\\\notag
&\lesod\nnodjt{\eta I\pludp\jsummup\li_n\AnAnz}^{\utfrac1{2q}\minu\frac12}\negc
\nnodjt{\zi I\pludp\jsummup\BnzBn}^{\utfrac1{2r}\minu\frac12}\negc
\nnodjd{\jsummup\li_n\AnAnz}^{\utfrac12\minu\frac1{2q}}\!
\nnoduu{\jsummup\!\BnzBn}^{\utfrac12\minu\frac1{2r}}\\
&\,\,\,\x\nnodut{\zzmddt{\jsummup\!\AnzAn}^\uvifrac1{2q}\negp X\! \zzmddc{\jsummup\!w_n\BnBnz}^\uvifrac1{2r}}_s\lessd
\nnodut{\zzmddt{\jsummup\!\AnzAn}^\uvifrac1{2q}\negp X\! \zzmddc{\jsummup\!w_n\BnBnz}^\uvifrac1{2r}}_s,\label{finis3}
\end{align}\vsplinc
where \ineqreff{finis1} follows from property $\nnnoo{CYD}_s\lesu\nnnoo{C}\nnnoo{Y}_s\nnnoo{D}$ for $C,D\inu\BH$ and $Y\inu\ccs,$ while \eqreff{finis2} is based on the property $\nnnoo{A^p}=\nnnoo{A}^p$ for positive operator $A$ and $p>0.$ The first inequality in \eqref{finis3} is based on \ineqreff{weghtedCSforCpAlt}, while the last inequality in \eqref{finis3} follows from the fact that $\nnnojjd{\eta I\pludp\TSjsummup\li_n\AnAnz}^{\uufrac1{2q}\minu\frac12}\! \nnnojtd{\TSjsummup\li_n\AnAnz}^{\uufrac12\minu\frac1{2q}}\lescj
\nnnojjd{\eta I\pludp\TSjsummup\li_n\AnAnz}^{\uufrac1{2q}\minu\frac12}\!
\nnnojjd{\eta I\pludp\TSjsummup\li_n\AnAnz}^{\uufrac12\minu\frac1{2q}}\eqpj1,$ and similarly
$\nnnojdt{\pozd\zi I\pludp\TSjsummup\BnzBn}^{\uufrac1{2r}\minu\frac12}\! \nnnojuu{\TSjsummup\BnzBn}^{\uufrac12\minu\frac1{2r}}\lespj1.$ By taking $\sup_{\eta,\zi>0}$ in inequalities \eqref{finis1}-\eqref{finis3} we finally get \ineqreff{weghtedCSforCp} .
\eeproo

\bbr
\ineqreff{weghtedCSforCp} essentially extends Ineq.\ (5) in the proof of \cite[Th.~2.1]{J99},
(which is equivalent to Ineq.\ (5) in \cite[Th.~2.1]{J99}), as well as  Ineq.\ (38) in the proof
of \cite[Th.~3.3]{J05} (which is equivalent to the main Ineq.\ (24) in  \cite[Th~3.3]{J05}).
\eer

In certain cases the square sumability for all of the four operator sequences appearing in Theorem
\ref{CSinCp} do not need to be required. The algebraic form of the related norm inequalities is also simpler, so in the next corollary we reformulate them in a more appropriate form.

\bbc\label{novkorolar}
Let $\!\nizjd{\lambda},\!\nizjd{w}\!\!$ be sequences in $\!\Rpz\!$ and $X\inu\ccs.$ \nlin
{\rm(}\negj a\pozj{\rm)} Let also $1\lesu s\lesu 2$ and $\nizjd{A},\odu{w_n^\uujd\!\Bnuz}_\nji\inuu\ldsBH.$ \nlin
{\rm(}\negj a1\pozj{\rm)} If additionally $\odjd{\li_n^\utjd\!\Anuz}_\nji\inuu\ldsBH,$
then there exists $\negc\TSCssummuc\li_n^{1\mindj\frac1s}\!w_n^\frac12\!\AnXBn\!$ and
\vspbbald\notag
&\nnodcd{\Cssumm\!\!\li_n^{1\mindj\frac1s}\!w_n^\frac12\!\AnXBn}_{\negd s}
\lessd\nnodod{\jsummuc\li_n\AnAnz}^{\negd1\mindj\frac1s}\!
\nnodct{\zzmddt{\jsummup\AnzAn}^{\uvifrac1s\mindj\frac12}\negp X\!
\zzmddc{\jsummup w_n\BnBnz}^\uvifrac12}_{\negd s}\!.
\end{align}\vsplinc
{\rm(}\negj a2\pozj{\rm)} If alternatively $\nizjd{B}\inuu\ldsBH,$ then there exists
$\negc\TSCssummuc\!w_n^{\ufrac1s\mindj\frac12}\!\AnXBn\!$ and
\vspbbald\notag
&\nnodcd{\Cssumm\!\!w_n^{\ufrac1s\mindj\frac12}\!\AnXBn}_{\negd s}
\lessd \nnodjd{\jsummup\BnzBn}^{\negd1\mindj\frac1s}\!
\nnodct{\zzmddt{\jsummup\AnzAn}^\uvifrac12\negp X\!
\zzmddc{\jsummup w_n\BnBnz}^{\uvifrac1s\mindj\frac12}}_{\negd s}\!.
\end{align}\vsplinc
{\rm(}\negj b\pozj{\rm)} Similarly, let us consider the case $s\gesu 2$ and $\odjd{\li_n^\utjd\!\Anuz}_\nji,\nizjd{B}\inuu\ldsBH.$ \nlin
{\rm(}\negj b1\pozj{\rm)} If also $\nizjd{A}\inuu\ldsBH,$  then there exists $\negc\TSCssummuc\!\li_n^{\ufrac12\mindj\frac1s}\!\AnXBn\!$ and
\vspbbald\notag
&\nnodcd{\Cssumm\!\!\li_n^{\ufrac12\mindj\frac1s}\!\AnXBn}_{\negd s}
\lessd\nnodod{\jsummuc\li_n\AnAnz}^{\utfrac12\mindj\frac1s}\!
\nnodjd{\jsummup\BnzBn}^{\utfrac12}\!
\nnodcd{\zzmddt{\jsummup\AnzAn}^{\uvifrac1s}\!X}_{\negd s}\!.
\end{align}\vsplind
{\rm(}\negj b2\pozj{\rm)} If alternatively $\odu{w_n^\uujd\!\Bnuz}_\nji\inuu\ldsBH,$
then there exists $\negc\TSCssummuc\!\li_n^\frac12 w_n^{\ufrac1{s}}\negj\AnXBn\!$ and
\vspbbald\notag
&\nnodcd{\Cssumm\!\!\li_n^\frac12 w_n^\ufrac1s\negj\AnXBn}_{\negd s}
\lessd \nnodod{\jsummuc\li_n\AnAnz}^\utfrac12\negd
\nnodjd{\jsummup\BnzBn}^{\utfrac12\mindj\frac1s}\!
\nnoddt{X\!\zzmddc{\jsummup w_n\BnBnz}^\uvifrac1s}_{\negd s}\!.
\end{align}\vsplinx
\eec
\bbproo
The inequality in (a1) (resp.\ (a2)) represents  the special case
$q\dod \frac{s}{2-s},$ $r\dod 1$ (resp.\ $q\dod 1,$ $r\dod\frac{s}{2-s}$) of \ineqreff{weghtedCSforCpAlt},
while inequality in (b1) (resp.\ (b2)) represents the special case
$q\dod \frac{s}2,$ $r\dod \iii$ (resp.\ $q\dod \iii,$ $r\dod\frac{s}2$)
of the same \ineqreff{weghtedCSforCpAlt}.
\eeproo

\bbc\label{korolar}
 Let $q,r,s\gesuu1\!$ satisfy $\!\frac1{2q}\plu\frac1{2r}\equt\frac1s\!$ and $X\inuu\ccs.$
If  $\!\nizjd{\lambda},\!\nizjd{w}\!\!$ are sequences in $\!\Rpz$ and $\nizjd{\li_n^{\negc\odu{1\minu q}\negd/\negd\odu{2q}}\negc A},\odjd{\li_n^{\negc1\negd/\negd\odu{2q}}\negc\Anuz}_\nji, \nizjd{w_n^{\mindt1\negd/\negd\odu{2r}}\!B},\odjd{w_n^{\negd\odu{r\minu1}\negd/\negd\odu{2r}}\!\Bnuz}_\nji\!
\inuu\ldsBH$  then there exists  $\TSCssummuc\AnXBn$ and
\vspbbald\notag
&\sup_{\eta,\zi\negd\IIS>\negj0}\negj
\nnoduc{\zzmdct{\eta I\plutp\jsummup\li_n^{\utfrac1q}\AnAnz}^{\uvifrac1{2q}\minu\frac12}\negp
\Cssumm\!\!A_n^\fff X\negd B_n^\fff\!
\zzmduc{\zi I\plutp\jsummup w_n^{\mintj\frac1r}\!\BnzBn}^{\uvifrac1{2r}\minu\frac12}}_{\negd s}\\
&\lesod\nnodct{\zzmddc{\jsummup\li_n^{\utfrac1q\minu1}\!\!\AnzAn}^\uvifrac1{2q}\negp X\!
\zzmddc{\jsummup w_n^{\negc1\minu\frac1r}\negc\BnBnz}^\uvifrac1{2r}}_{\negd s}\negd,\label{altotag}\\
&\nnoduu{\negd\Cssummuc\AnXBn}_s
\lessd\nnodod{\jsummuc\li_n^{\utfrac1q}\AnAnz}^{\utfrac12\minu\frac1{2q}}\!
\nnodjd{\jsummup w_n^{\mintj\frac1r}\!\BnzBn}^{\utfrac12\minu\frac1{2r}}\!
\nnodct{\zzmddt{\jsummup\li_n^{\utfrac1q\minu1}\negp\AnzAn}^\uvifrac1{2q}\negp X\!
\zzmddc{\jsummup w_n^{\negc1\minu\frac1r}\!\BnBnz}^\uvifrac1{2r}}_s\!.\label{altotagg}
\end{align}\vsplindp
If  $\!\nizjd{\gi},\!\nizjd{\rho}\!\!$ are sequences in $\!\Rpz$ and $\nizjd{A},\odjd{\gi_n^\utfrac{q}{q-1}\negc\Anuz}_\nji, \nizjd{B},\odjd{\rho_n^r\negj\Bnuz}_\nji\!
\inuu\ldsBH$  then there exists  $\TSCssummuc\gi_n\rho_n\AnXBn$ and
\vspbbald\label{weghtedCSforCpsmena}
&\sup_{\eta,\zi\negd\IIS>\negj0}\negj
\nnodct{\zzmdct{\eta I\plutp\jsummup \gi_n^\utfrac{2q}{q-1}\negc\AnAnz}^{\uvifrac1{2q}\minu\frac12}\negp
\Cssumm\!\!\gi_n \rho_n\AnXBn\!
\zzmduc{\zi I\plutp\jsummup\BnzBn}^{\uvifrac1{2r}\minu\frac12}}_{\negd s}
\lesst\nnodct{\zzmddc{\jsummup\AnzAn}^\uvifrac1{2q}\negp X\!\zzmddc{\jsummup \rho_n^{2r}\BnBnz}^\uvifrac1{2r}}_{\negd s}\negc,\\
&\nnodcd{\Cssumm\!\gi_n \rho_n\AnXBn}_{\negd s}
\lessd\nnodod{\jsummuc\gi_n^\utfrac{2q}{q-1}\negc\AnAnz}^{\uufrac12\minu\frac1{2q}}\!
\nnodjd{\jsummup\BnzBn}^{\uufrac12\minu\frac1{2r}}\!
\nnodct{\zzmddt{\jsummup\AnzAn}^\uvifrac1{2q}\negp X\!
\zzmddc{\jsummup \rho_n^{2r}\BnBnz}^\uvifrac1{2r}}_{\negd s}\!.\label{weghtedCSforCpAltsmena}
\end{align}\vsplinc
\eec\vspc
\bbproo
The proof for \ineqreff{altotag} (resp.\ \ineqreff{altotagg})
follows directly from \ineqreff{weghtedCSforCp} (resp.\ \ineqreff{weghtedCSforCpAlt}) if we take $\odju{\li_n^{\uufrac1{2q}\minu\frac12}\negp A_n, w_n^{\mintj\frac1{2r}}\!B_n}$ instead of
$\odu{A_n, B_n}$ for any $n\inN.$

Similarly, \ineqreff{weghtedCSforCpsmena} (resp.\ \ineqreff{weghtedCSforCpAltsmena}) obtains
directly from \ineqreff{weghtedCSforCp} (resp.\ \ineqreff{weghtedCSforCpAlt}) if we take
$\li_n\dodu\gi_n^\utfrac{2q}{q-1}$ and $w_n\dodu \rho_n^{2r}$ for any $n\inN.$
\eeproo

\section{Applications to some asymmetrically square summable families of operators}

The next two theorems presents the first applications of Cauchy--Schwarz norm inequalities
in Schatten--von Neumann ideals for \awss sequences of operators.

\bbt
Let  $q,r,s\gesu1$ satisfy $\frac1{2q}\plu\frac1{2r}=\frac1{s},$ $\nizjd{e}\!$ and $\nizjd{f}\!$ are \ONB in $\HH$ and $X\inu\ccs.$
 If $\nizjd{\li}$ and $\nizjd{\rho}$ are summable sequences in $\Rpz,$ then
\vspbbald\label{EkstremWeghtedCSforCp}
&\sup_{\eta\negd\IIS>\negj0}
\nnodct{\zzmduc{\eta I\pludt\zzmdjc{\summut\li_n}e_1\zoxtd e_1}^{\uvifrac1{2q}\minu\frac12}\!
\Cssumm\negc\li_n^{\ucfrac12\mindj\frac1{2q}}\!\rho_n^\uufrac1{2r}\!\zzpou{Xf_1,e_n}f_n\zoxtd e_1}_s
\lesst\zzmddp{\summud\rho_n}^\uvifrac1{2r}\!\!\nnnou{X\negj f_1}\!,\\
&\nnodct{\Cssumm\negc\li_n^{\ucfrac12\mindj\frac1{2q}}\!\rho_n^\frac1{2r}
\!\zzpou{Xf_1,e_n}f_n\zoxtd e_1}_s
\lessd\zzmddc{\summut\li_n}^{\uvifrac12\minu\frac1{2q}}\negp
\zzmddc{\summud\rho_n}^\uvifrac1{2r}\!\!\nnnou{X\negj f_1}\!.\label{EkstremWeghtedCSforCpAlt}
\end{align}\vsplind
Consequently, if $\nizjd{\gi}\inu\ell^{\ufrac{2q}{q-1}}_{\IIS\eN}\!\!$  and
$\nizjd{w}\inu\ell^{2r}_{\IIS\eN}$ for some sequences $\nizjd{\gi},\nizjd{w}\!$ in $\Rpz,$ then
\vspbbald\label{EkstremWeghtedCSforCpAlt1}
&\sup_{\eta\negd\IIS>\negj0}
\nnodct{\zzmduc{\eta I\pludt\zzmdjc{\summut \gi_n^\utfrac{2q}{q-1}}e_1\zoxtd e_1}^{\uvifrac1{2q}\minu\frac12}\negp
\Cssumm\negc\gi_n w_n\!\zzpou{Xf_1,e_n}f_n\zoxtd e_1}_s
\lesst\zzmddc{\summut w_n^{2r}}^\uvifrac1{2r}\!\!\nnnou{X\negj f_1}\!,\\
&\nnodct{\Cssumm\negc\gi_n w_n
\!\zzpou{Xf_1,e_n}f_n\zoxtd e_1}_s
\lessd\zzmdjc{\summut\gi_n^\utfrac{2q}{q-1}}^{\usfrac{q-1}{2q}}\negp
\zzmddc{\summut w_n^{2r}}^\uvifrac1{2r}\!\!\nnnou{X\negj f_1}\!,\label{EkstremWeghtedCSforCpAlt2}
\end{align}\vsplinx
\eet
\bbproo
As $\odjjd{e_n\zoxtd e_1}_\nji\inu\ldsBH\captj\ldstliBH$ and $\odjjd{f_n\zoxtd f_1}_\nji\inu \ldsBH\captj\ldstroBH$ based on Example \ref{osnovniprimer},
then the application of \ineqreff{weghtedCSforCp} to sequences $\nizjd{A}$ and
$\nizjd{B},$ given by $A_n\dodu e_n\zoxtd e_1,B_n\dodu f_n\zoxtd f_1$ for any $n\inN,$ gives
\vspbbald\notag
&\negxx\sup_{\eta\negd\IIS>\negj0}
\nnodct{\zzmduc{\eta I\pludt\zzmdjc{\summut\li_n}e_1\zoxtd e_1}^{\uvifrac1{2q}\minu\frac12}\negp
\Cssumm\negc\li_n^{\ucfrac12\mindj\frac1{2q}}\!\rho_n^\uufrac1{2r}\!\zzpou{Xf_1,e_n}f_n\zoxtd e_1}_s\\
&\negxx=\sup_{\eta,\zi\negd\IIS>\negj0}
\nnodct{\zzmduc{\eta I\pludt\zzmdjc{\summut\li_n}e_1\zoxtd e_1}^{\uvifrac1{2q}\minu\frac12}\negp
\Cssumm\negc\li_n^{\ucfrac12\mindj\frac1{2q}}\!\rho_n^\uufrac1{2r}\!
\odjjd{e_n\zoxtd e_1}\negj X\!\odjjd{f_n\zoxtd f_1}\!\odu{1\plu\zi}^{\utfrac1{2r}\minu\frac12}}_s
\lesst\nnoddt{X\!\zzmdpc{\zzmdjp{\summud\rho_n}f_1\zoxtd  f_1}^\uvifrac1{2r}}_s\!. \label{EkstremWeghtedCSforCpDok}
\end{align}\vsplinc
Here we rely on \ineqreff{ektremPrimer} and \ineqreff{ektremPrimer2} to see
that $\TSjsummup\li_n\AnAnz=\zzmjod{\sumN\li_n} e_1\zoxtj  e_1,$ $\TSjsummup\AnzAn=\TSjsummup\BnzBn=I$
and $\TSjsummup \rho_n\BnBnz=\zzmjod{\sumN \rho_n} f_1\zoxtd  f_1.$
Also, \eqreff{EkstremWeghtedCSforCpDok} is due to $\sup_{\zi\negd\IIS>\negj0}\odu{1\plu\zi}^{\utfrac1{2r}\minu\frac12}
=\lim_{\zi\to0^+}\odu{1\plu\zi}^{\utfrac1{2r}\minu\frac12}=1.$ So, to obtain \ineqreff{EkstremWeghtedCSforCp} from \ineqreff{EkstremWeghtedCSforCpDok}, we just have to take into account that for the rank one ortoprojector $f_1\zoxtd f_1$ we have
$\nnnoou{X\zzmojd{f_1\zoxtd  f_1}^{\uufrac1{2r}}}_s
=\nnnoou{X\zzmojd{f_1\zoxtd  f_1}}_s= \nnnou{f_1\zoxtd X\negj f_1}_s\negd=\nnnou{X\negj f_1}.$

Similarly, the proof of \ineqreff{EkstremWeghtedCSforCpAlt} is based on
\vspbbalzd
\nnodct{\Cssumm\negc\li_n^{\ucfrac12\mindj\frac1{2q}}\!\rho_n^\frac1{2r}\!
\zzpou{Xf_1,e_n}f_n\zoxtd e_1}_s
\eqpd\nnodct{\Cssumm\negc\li_n^{\ucfrac12\mindj\frac1{2q}}\!\rho_n^\frac1{2r}\!
\odjjd{e_n\zoxtd e_1}\negj X\!\odjjd{f_n\zoxtd f_1}}_s
\lessd\zzmdjc{\summut\li_n}^{\uvifrac12\minu\frac1{2q}}\negp
\zzmddp{\summud\rho_n}^\uvifrac1{2r}\!\!\nnnou{X\negj f_1}\!,
\end{align*}\vsplinc
where we applied \ineqreff{weghtedCSforCpAlt} to the same families $\nizjd{A}\!$ and $\nizjd{B}\!$ appearing in the proof of \ineqreff{EkstremWeghtedCSforCp}.

Finally, \ineqreff{EkstremWeghtedCSforCpAlt1} (resp.\ \ineqreff{EkstremWeghtedCSforCpAlt2}) follows directly from
\ineqreff{EkstremWeghtedCSforCp} (resp.\
\ineqreff{EkstremWeghtedCSforCpAlt}) by taking $\li_n\dodu \gi_n^\utfrac{2q}{q-1}$
and $\rho_n\dodu w_n^{2r}$ for any $n\inN.$
\eeproo

\smallskip

We conclude this paper with application of \ineqreff{weghtedCSforCp} and \ineqreff{weghtedCSforCpAlt}
to sequences of operators related to hypercontractive operators.

\bbt
If $C$ (resp.\ $D$) is $\N$-cohypercontractive (resp.\ $\M$-cohypercontractive) for some $\N,\M\inu\eN,$ then the operator sequence $C_{\!\NN}\dodu\odjte{\TSsqrt{\GCzCI{}\!}\,C^n\!\TSsqrt{\GCCzI{\NN}\!}}_\uunoi\!$ is $\li^{\negd*}\negd$-\awss for $\li\dodu\zzvjuc{\tbinomu{n\plu\NN\minuu1}{\NN\minuu1}}_\uunoi,$
$D_{\!\MM}\dodu\odjte{\TSsqrt{\GDzDI{}\!}\,D^n\!\TSsqrt{\GDDzI{\MM}\!}}_\uunoi\!$ is $\rho^*\negd$-\awss for $\rho\dodu\zzvjuc{\tbinomu{n\plu\MM\minuu1}{\MM\minuu1}}_\uunoi,$
i.e.\ $C_{\!\NN}\inu\ldsBH\cap\ldstliBH$ and $D_{\!\MM}\inu\ldsBH\cap\ldstroBH,$ and if  $q,r,s\gesu1$ satisfy $\frac1{2q}\plu\frac1{2r}=\frac1{s},$ then
\vspbbald\notag
&\!\sup_{\eta,\zi\IIS>\negj0}\LLndc{\odjdt{\eta I\pluu\GCzCI{}}^{\ucfrac1{2q}\minu\frac12}\negc
\Cssummoup\tbinomu{n\plu\NN\minuu1}{\NN\minuu1}^{\ucfrac12\mindj\frac1{2q}}\!
\tbinomu{n\plu\MM\minuu1}{\MM\minuu1}^\ucfrac1{2r}\negp
\TSsqrt{\GCzCI{}\!}\,C^n\!\TSsqrt{\GCCzI{\NN}\!}\,X\negd\TSsqrt{\GDzDI{}\!}\,D^n\!\TSsqrt{\GDDzI{\MM}\!}}\\
&\quad\pozs\x\RRndt{\odjdt{\zi I\pluu\GDDzI{\MM}}^{\ucfrac1{2r}\minu\frac12}}_s
\lespt\nnnojut{\odjjt{\GCCzI{\NN}}^\ucfrac1{2q}\negc X\!\odjjt{\GDzDI{}}^\ucfrac1{2r}}_s\negj, \label{ExWeghtedCSforCp}\\
&\nnodcj{\Cssummoup\tbinomu{n\plu\NN\minuu1}{\NN\minuu1}^{\ucfrac12\mindj\frac1{2q}}\!
\tbinomu{n\plu\MM\minuu1}{\MM\minuu1}^\ucfrac1{2r}\negp
\TSsqrt{\GCzCI{}\!}\,C^n\!\TSsqrt{\GCCzI{\NN}\!}\,X\negd\TSsqrt{\GDzDI{}\!}\,D^n\!\TSsqrt{\GDDzI{\MM}\!}\,}_s\notag\\
&\lesu\nnnoju{\GCzCI{}}^{\uufrac12\minu\frac1{2q}}\!\nnnoju{\GDDzI{\MM}}^{\uufrac12\minu\frac1{2r}}\!
\nnnojut{\odjjt{\GCCzI{\NN}}^\ucfrac1{2q}\negc X\!\odjjt{\GDzDI{}}^\ucfrac1{2r}}_s\negd. \label{ExWeghtedCSforCpAlt}
\end{align}\vsplinxi
\eet
\bbproo 
The proof on \ineqreff{ExWeghtedCSforCp} is based on the following estimates:
\vspbbald\notag
&\sup_{\eta,\zi\negd\IIS>\negj0}\LLndc{\odjdt{\eta I\pluu\GCzCI{}}^{\ucfrac1{2q}\minu\frac12}\negc
\Cssummoup\tbinomu{n\plu\NN\minuu1}{\NN\minuu1}^{\ucfrac12\mindj\frac1{2q}}\!
\tbinomu{n\plu\MM\minuu1}{\MM\minuu1}^\ucfrac1{2r}\negp
\TSsqrt{\GCzCI{}\!}\,C^n\!\TSsqrt{\GCCzI{\NN}\!}\,X\negd\TSsqrt{\GDzDI{}\!}\,D^n\!\TSsqrt{\GDDzI{\MM}\!}}\\
&\quad\pozs\x\RRndt{\odjdt{\zi I\pluu\GDDzI{\MM}}^{\ucfrac1{2r}\minu\frac12}}_s
\lesup\sup_{\eta,\zi\negd\IIS>\negj0}\LLnd{\negc\zzmdtj{\eta I\pludp\jsummo\!\!\tbinomu{n\plu\NN\minuu1}{\NN\minu1} \negd\TSsqrt{\GCzCI{}\!}\,C^n\GCCzI{\NN}C^{*n}\! \TSsqrt{\GCzCI{}\!}}^{\uvifrac1{2q}\minu\frac12}}\label{ExWeghtedCSforCpDok1}\\
&\quad\pozs\x\Cssummoup\tbinomu{n\plu\NN\minuu1}{\NN\minuu1}^{\ucfrac12\mindj\frac1{2q}}\!
\tbinomu{n\plu\MM\minuu1}{\MM\minuu1}^\ucfrac1{2r}\negp
\TSsqrt{\GCzCI{}\!}\,C^n\!\TSsqrt{\GCCzI{\NN}\!}\,X\negd\TSsqrt{\GDzDI{}\!}\,D^n\!\TSsqrt{\GDDzI{\MM}\!}\notag\\
&\quad\pozs\x\RRnd{\zzmdcj{\zi I\pludp\jsummo\!\! \TSsqrt{\GDDzI{\MM}\!}\,D^{\negj*n}\GDzDI{}D^n\!\TSsqrt{\GDDzI{\MM}\!}}^{\uvifrac1{2r}\minu\frac12}\!}_s\notag\\
&\lesot\nnodct{\zzmdtj{\jsummo\! \negd\TSsqrt{\GCCzI{\NN}\!}\,C^{*n}\GCzCI{}C^n\!\TSsqrt{\GCCzI{\NN}\!}}^\ucfrac1{2q}\negc X
\negd\zzmdtj{\jsummo\!\!\tbinomu{n\plu\MM\minuu1}{\MM\minu1} \negd\TSsqrt{\GDzDI{}\!}\,D^n\GDDzI{\MM}D^{*n}\!\TSsqrt{\GDzDI{}\!}}^{\uvifrac1{2r}}}_s\negj \label{ExWeghtedCSforCpDok2}\\
&\lesot\nnnojut{\odjjt{\GCCzI{\NN}}^\ucfrac1{2q}\negc X\!\odjjt{\GDzDI{}}^\ucfrac1{2r}}_s\negj. \label{ExWeghtedCSforCpDok3}
\end{align}\vsplinc
Here, \ineqreff{ExWeghtedCSforCpDok1} is due to estimates \eqref{hiperkontr2} and
\eqref{hiperkontr3}, by using also the fact that the function $t\mapsto -t^\ai$ is operator monotone on $\eRpo$ if $\ai\inu[-1,0),$ for $\ai\dodu\frac1q\mindj1$ and $\ai\dodu\frac1r\mindj1,$ all combined with the double monotonicity property \eqref{duplamonotonost} for $\nnnou{\cdot}_s$ norm.
\ineqreff{ExWeghtedCSforCpDok2} follows from \ineqreff{weghtedCSforCp} applied to $\li_n\dod \tbinomu{n\plu\NN\minuu1}{\NN\minuu1},$ $w_n\dodu \tbinomu{n\plu\MM\minuu1}{\MM\minuu1}$ (i.e.\ $w\dod \rho$) and to $(\TSsqrt{\GCzCI{}\!}\,C^n\!\TSsqrt{\GCCzI{\NN}\!}\,,
X,\negd\TSsqrt{\GDzDI{}\!}\,D^n\!\TSsqrt{\GDDzI{\MM}})$ instead of $(A_n,X,B_n)$ therein, while \ineqreff{ExWeghtedCSforCpDok3} is based on inequalities \eqref{hiperkontr1} and \eqref{hiperkontr4},
the fact that the function $t\mapsto t^\ai$ is operator monotone on $[0,\iip)$ if $\ai\inu(0,1]$
(see \cite[Ex.~2.5.9(3)]{Hi}) for $\ai\dodu\frac1q$ and $\ai\dodu\frac1r,$ combined with the double monotonicity property \eqref{duplamonotonost}.

Similarly, to prove \ineqreff{ExWeghtedCSforCpAlt} we just apply \ineqreff{weghtedCSforCpAlt} instead of \ineqreff{weghtedCSforCp} to $(\negj\TSsqrt{\GCzCI{}\!}\,C^n\!\TSsqrt{\GCCzI{\NN}\!}\,,
X,\allowbreak\negd\TSsqrt{\GDzDI{}\!}\,D^n\!\TSsqrt{\GDDzI{\MM}})$ (instead of $(A_n,X,B_n)\negd$)
for the same $\li,\rho$  as in the proof of \eqref{ExWeghtedCSforCp},
by repeating the same arguments.
\eeproo

\bbr Note that \ineqreff{EkstremWeghtedCSforCp} (resp.\ \ineqreff{EkstremWeghtedCSforCpAlt})
provides the example for the application of \ineqreff{altotag} (resp.\ \ineqreff{altotagg}) in the
Corollary \ref{korolar} to sequences $(A_n,B_n)\dodu (\li_n^{\ucfrac12\mindj\frac1{2q}}\!e_n\zoxtd e_1,\rho_n^\uufrac1{2r}\!f_n\zoxtd f_1).$

In the same spirit, \ineqreff{ExWeghtedCSforCp} (resp.\ \ineqreff{ExWeghtedCSforCpAlt})
also provides the example for the application of \ineqreff{altotag} (resp.\ \ineqreff{altotagg}) in
Corollary \ref{korolar} to the sequences
$\odu{A_n,\negd B_n,\negd \li_n,\negd w_n}\dodu \bigl(\!\tbinomu{n\plu\NN\minuu1}{\NN\minuu1}^{\ucfrac12\mindj\frac1{2q}}\negc \TSsqrt{\GCzCI{}\!}\,C^n\!\TSsqrt{\GCCzI{\NN}\!}\,,\allowbreak
\tbinomu{n\plu\MM\minuu1}{\MM\minuu1}^\uufrac1{2r}\negc\TSsqrt{\GDzDI{}\!}\,D^n\!\TSsqrt{\GDDzI{\MM}\!}\,,
\negd\tbinomu{n\plu\NN\minuu1}{\NN\minuu1},\negd\tbinomu{n\plu\MM\minuu1}{\MM\minuu1}\!\bigr).$
\eer

\section*{Funding} 
The research of the authors was partially funded by Faculty of Mathematics University of Belgrade through the grant by the Ministry of Science, Technological Development and Innovation of the Republic of Serbia (the contract no 451-03-136/2025-03/200104)

\section*{Data availability statement} No datasets were generated or analyzed during the current study.

\section*{Declarations}
\subsection*{Competing interests} The authors declare no competing interests.
\\

\textbf{ORCID}\qquad Danko R.\ Joci\'c \qquad  http://orcid.org/0000-0003-2084-7180

\qquad\qquad\quad\pozxx\pozt Mihailo  Krsti\'c \qquad http://orcid.org/0000-0003-3575-3216

\qquad\qquad\quad\pozx  Milan Lazarevi\'c \qquad  http://orcid.org/0000-0003-1408-5626

\qquad\quad\pozx\pozj\,\,\, Stevan Mila\v sinovi\'c \qquad http://orcid.org/0000-0003-3088-0902


\end{document}